\documentclass[12pt]{amsart}
\usepackage{amsmath,amssymb,amsthm,epsfig}
\def\torus{\mathbb T}
\def\real{\mathbb R}
\def\complex{\mathbb C}
\def\integer{\mathbb Z}
\def\const{\text{const}}

\def\Iso{\qopname\relax o{Iso}}
\newtheorem{theorem}{Theorem}
\newtheorem{lemma}{Lemma}
\theoremstyle{remark}
\newtheorem{ex}{Example}

\newtheorem{conjecture}{Conjecture}
\newtheorem{remark}{Remark}
\author{I.\,Dynnikov}
\author{ S.\,Novikov}
\address{Dept.\ of Mech.\ and Math., Moscow State University,
Moscow 119992 GSP-2, Russia}
\email{dynnikov@mech.math.msu.su}
\address{Landau Institute for theoretical physics, Kosygina str.~2,
Moscow 119334, Russia; IPST, University of Maryland, College Park,
MD 20742, USA}
\email{novikov@glue.umd.edu}
\title{Topology of quasiperiodic functions on the plane}
\thanks{The work of I.\,Dynnikov was supported in part by
Russian Foundation for Basic Research (grant no.~02-01-00659);
the work of S.\,Novikov was supported in part by
the Council of the Russian Academy of Science (grant ``Mathematical
methods of nonlinear dynamics'')}

\begin{document}

\maketitle

\begin{abstract}This article describes a topological theory of
quasiperiodic functions on the plane. The development of
this theory was started (in different
terminology) by the Moscow topology group in early 1980s.
It was motivated by the needs of solid state physics, as a partial
(nongeneric) case of Hamiltonian foliations of  Fermi
surfaces with multivalued Hamiltonian function~\cite{N}. The
unexpected discoveries of their topological properties that were made in
1980s \cite{Z,N1} and 1990s \cite{D,D',D1} have finally led to
nontrivial physical conclusions \cite{NM,NM1} along the lines of
the so-called geometric strong magnetic field limit
\cite{LAK}. A very fruitful new point of view comes from
the reformulation of that problem
in terms of quasiperiodic functions and an extension to higher dimensions
made in 1999 \cite{N2}. One may say that,
for single crystal normal metals put in a magnetic field,
the semiclassical trajectories of electrons in the space of
quasimomenta are exactly the level lines of the quasiperiodic function with
three quasiperiods that is the dispersion relation restricted to a
plane orthogonal to the magnetic field. General studies of the
topological properties of levels of quasiperiodic functions on the
plane with any number of quasiperiods were started in 1999 when
certain ideas were formulated for the case of four quasiperiods
\cite{N2}. The last section of this work contains a complete
proof of these results based on the technique developed in \cite{D2,D3}.
Some new physical applications of the general problem were found
recently \cite{M}.
\end{abstract}

\section{Quasiperiodic functions}

Let $\torus^n=\real^n/\integer^n$ denote the
$n$-dimensional torus,
$\nu\!:\real^n\rightarrow\torus^n=\real^n/\integer^n$ the standard
projection.

We say that a real smooth function $\varphi(y)=\varphi(y^1,\dots,y^k)$
on the $k$-plane $\real^k$ is \emph{quasiperiodic} with $n$
quasiperiods (frequencies) if it can be represented in the form
$\varphi(y)=f(x(y))$:
\begin{equation}
\varphi=f\circ\nu\circ \iota,
\end{equation}
where $\iota:\real^k\rightarrow\real^n$ is an affine imbedding:
$$x^s=a^s_ry^r+x^s_0,$$
$f=f(x):\torus^n\rightarrow\real$ is some
smooth function, and $n\geqslant 2$ is the minimal possible integer for
which such a function $f$ and an affine imbedding $\iota$ exist. Here
$s=1,\dots,n$ and $r=1,\dots,k$. In the theory of quasicrystals, people call
the space $\real^k$ (where $k=2,3$) the
\emph{physical space} and the space $\real^n$
the \emph{superspace}. Every $n$-periodic function $f(x)$ generates a
family of \emph{descendants}, which are obtained by varying
the initial vector
$x_0=(x^1_0,\dots,x^n_0)$ in the superspace $\real^n$.
Any two descendants $\varphi_1(y)$, $\varphi_2(y)$
of the same function $f$ are said to be
\emph{related}. They have the same frequencies
and obtain the following property: for any $\varepsilon>0$,
there is a shift $y\mapsto y+a$ in the physical space
such that the shifted function
$\varphi_2(y+a)$ is $\varepsilon$-close to $\varphi_1(y)$:
$$|\varphi_2(y+a)-\varphi_1(y)|<\varepsilon\qquad\forall y\in\real^k.$$

 Any
linear function $\lambda:\real^k\rightarrow\real$ of the form
$\lambda(y)=\ell(x(y))$ or $\lambda=\ell\circ \iota$, where the
linear function $\ell:\real^n\rightarrow\real$ belongs to the dual
(or ``reciprocal'') lattice $(\integer^n)^*$ (\emph{i.e.}, we have
$\ell(\integer^n)\subset\integer$) is called a \emph{frequency} of
$\varphi$. The set of all frequencies form a free abelian group
with $n$ natural  generators
$\lambda_1=\ell^1\circ\iota,\dots,\lambda_1=\ell^n\circ\iota$ where
the functions
$\ell^s(x)=x^s$, $s=1,\ldots,n$, are dual to
the basic periods. We call this group
the \emph{group of frequencies}. It is a free
abelian subgroup $\Gamma^*$ of the dual vector space $\real^*$,
and it is the
same for the whole family of related quasiperiodic functions (descendants of the
same $n$-periodic function $f$).

Analytically, any $n$-periodic function can be presented in the form
of a trigonometric series

$$f(x)=\sum_{\ell\in(\integer^n)^*} c_\ell\exp\big(2\pi
i\ell(x)\big)$$

Therefore, any  quasiperiodic function can be presented in a
similar form:
$$\varphi(y)=\sum_{\lambda\in\Gamma^*}b_{\lambda}\exp\big(2\pi
i\lambda(y)\big)= \sum_mb_m\exp\Big(2\pi i
\sum_{s=1}^nm_s\lambda_s(y)\Big),$$
where $m=(m_1,\dots,m_n)
\in\integer^n$. By definition, the set of basic frequencies
$\lambda_s$ generates the space $\real^k$ over the field $\real$.

For the space $\real^k$ endowed with a Euclidean metric there is a natural
identification $\real^k\cong(\real^k)^*$, so the subgroup of
frequencies can be treated as a subgroup $\Gamma\cong\Gamma^*\subset
\real^k$ in the physical space $\real^k$.

There is an affine symmetry semigroup associated with each family
of related quasiperiodic functions. By definition, this semigroup
$\widetilde{G}$ consists of all affine transformations
$$g:\real^k\rightarrow \real^k$$
of the physical space
such that
$$g(\Gamma^*)\subset\Gamma^*,$$
where $\Gamma^*$ is treated as
a subset of the group of translations: $\Gamma^*\subset \real^k$.
For the Euclidean space $\real^k$ we define the \emph{symmetry group}
$G\subset\widetilde{G}$ consisting only of isometries $g$ such
that $g(\Gamma^*)=\Gamma^*$.

This group satisfies the general definition of a
quasicrystallographic group introduced by S.\,Novikov and
A.\,Veselov in 1980s in order to answer the question: what is the
symmetry of quasicrystals (see~\cite{Pi})? According to that
definition the intersection $G\bigcap\real^k\subset\Iso(\real^k)$
of a quasicrystallographic group
with the subgroup of translations should be a finitely generated free
abelian group. In our case it is exactly the group $\Gamma^*$. The
definition allows  the ``rotational'' quotient group
$G/(G\bigcap \real^k)\subset O_k$ to be infinite.
S.\,Piunikhin studied these
groups for $k=2,3$ in a series of works (see \cite{Pi}).

\begin{ex} Consider the two-dimensional case, $k=2$. Let $\theta$ be a
unimodular complex number, $|\theta |=1$, $\theta=\exp(i\psi)$,
satisfying the equation
$$P(\theta)=\theta^n+a_1\theta^{n-1}+\ldots+a_{n-1}\theta+1=0,$$
where all coefficients are  integer-valued, $a_s\in
\integer$, and we have $a_s=a_{n-s}$. The complex numbers
(or real two-vectors)
$\lambda_1=1, \lambda_2=\theta,\dots, \lambda_n=\theta^{n-1}$
generate a group of frequencies $\Gamma^*\subset
\complex=\real^2$  with nontrivial rotational symmetry $g\rightarrow
g\exp(i\psi)$. It is easy to find such a polynomial $P$ with
a root at $\theta=\exp(i\psi)$, where the
ratio $\psi/2\pi$ is irrational. There are very complicated
quasicrystallographic groups for $k=3$ (see \cite{Pi}).

\end{ex}

\section{Quasiperiodic functions in analysis, geometry, and
dynamical systems motivated by natural sciences}

\subsection{Quasiperiodic functions on the real line}

Consider the case $k=1$. In the XIX century,
one-dimensional quasiperiodic functions
with $n$ quasiperiods appeared in the theory of
completely integrable Hamiltonian systems of the classical
mechanics with $n$ degrees of freedom. According to so-called
Liouville's theorem, the integrability follows from the existence of
$n$ smooth independent pairwise commuting integrals of motion. If their
common level sets are compact, then the time dependence of the space
coordinates along a trajectory can be described by
quasiperiodic functions $x^r(t)$ with (at most) $n$ quasiperiods.
So, all studies of perturbations of
completely integrable systems should start with quasiperiodic
unperturbed background. A lot of fundamental work has been done in
this area (see \cite{DKN}).

\subsection{Quasiperiodic functions in the theory of nonlinear PDE}
Completely integrable PDE systems of the theory of solitons give rise
to quasiperiodic functions with $k>1$. There are very famous
($1+1$) PDE systems such as KdV ($u_t=6uu_x+u_{xxx}$) or sine-Gordon
($u_{tt}=u_{xx}+\sin \{u(x,t)\}$), which are completely integrable by the
so-called \emph{inverse scattering transform method for rapidly
decreasing initial values}. A countable number of continuous
families of exact
smooth real ``finite-gap'' solutions of these equations were discovered
in 1970s (see~\cite{DKN}). These solutions are quasiperiodic
functions in $x,t$, and depend on many parameters $a$, $a'$:
$$u(x,t)=F(xU+Vt+U_0;a)$$
for KdV, and
$$\exp\big(iu(x,t)\big)=F'(U'x+V't+U'_0;a')$$ for sine-Gordon. Here
$u(x,t)$ is real in both cases, $F,F'$ are $n$-periodic smooth
functions in $n$ variables (\emph{i.e.}, smooth functions on the real
$n$-torus). They can be expressed through special functions, namely,
theta-functions of a hyperelliptic Riemann surface
of genus~$n$. $U$, $U'$, $V$, $V'$ are the $n$-vectors of
periods of some Abelian differentials of the second kind
(see \cite{N1}). Let us
mention that, for the sine-Gordon system, the function $u=1/i\,\log
F'$ is generically a multivalued function on the ``real''
$n$-torus imbedded in the complex $2n$-dimensional Jacoby torus
associated with a complex hyperelliptic Riemann surface. Here we
have $k=2$. For famous completely integrable $(2+1)$ systems
(like KP, and others) one comes to quasiperiodic solutions of the form
$u(x,y,t)$, which are quasiperiodic functions
in $k=3$ physical variables. When studying the
dependence of the solution on so-called higher times one may arrive at any
value of $k$.

\subsection{Quasiperiodic functions and quasicrystals}

Completely different examples come from solid
state physics. In 1980s a new type of 2D and 3D media was
discovered. People named them ``quasicrystals''. The optical
analysis of the location of atoms gave an evidence for the group
of frequencies being incompatible with an ordinary crystal structure.
For example,
for $k=2$, the observed group of frequencies $\Gamma^*$ might be generated
by the 5th roots of unity:
$$\lambda_r=\eta^r\in\Gamma^*,\quad r=0,1,2,3,\quad\eta^5=1,\quad P(\eta)=0,$$
where
$$P(z)=z^4+z^3+z^2+z+1.$$

Recall that our extension of the idea of symmetry allows the
rotational symmetry to be even infinite.

There are two mathematical models of quasicrystals.
Let us think of atoms in
the physical space $\real^k$ as being located in a discrete set of
points $x_A$ such that
there exists a couple of positive ``radii'' $\rho_1$, $\rho_2$ with the properties:
\begin{itemize}
\item[a.]
We have $|x_A-x_{A'}|\geq\rho_2$ for all pairs $A,A'$ with $A\neq
A'$;
\item[b.]
For every point $x\in \real^k$, there exists a point $x_{A}$
such that $|x_A-x|\leqslant \rho_1$.
\end{itemize}
We call this set of points
\emph{quasiperiodic} if the distribution $\sum_A\delta(x-x_A)$ can be
decomposed into a Fourier series with finitely generated free
abelian group of frequencies $\Gamma^*$.

In another model, our
physical space $\real^k$ is endowed with
a ``quasiperiodic tiling''. This means the following:
\begin{itemize}
\item[a.] The space is
covered by countably many polytopes $P_B$,
$\real^k=\bigcup_B P_B$, where $P_B\bigcap P_{B'}$ is a face for any pair $B,B'$.
\item[b.]
Up to shift, there is only a finite number of different polytopes
$P_1,\dots,P_q$ among them.
\item[c.]
Let us associate some constant $c_q$ with every polytope $P_q$
and consider a function that is equal to $c_j$ everywhere
in the interior of any $P_B$ obtained from $P_j$ by a shift. We obtain a
piecewise constant function $c(x)$ in $x\in \real^k$ defined
(at a full measure set)
by our tiling and the choice of the constants $c_j$.
The tiling is said to be \emph{quasiperiodic} if, for every
choice of constants $c_j$, the function $c(x)$ is quasiperiodic,
\emph{i.e.}, can be presented in the form of a
trigonometric series with finitely generated free abelian group $\Gamma^*$ of
frequencies.
\end{itemize}

There is a famous tiling of
the plane $\real^2$ by rhombi of two types: one with
angles $\pi/5$ and $4\pi/5$,
and the other with angles $2\pi/5$ and $3\pi/5$.
It is called the Penrose tiling.
This tiling is quasiperiodic, which was discovered a few years later
after Penrose's original work
(see the history and details of the subject in~\cite{Pi}). An
interesting idea of ``local rules'' was developed by
physicists and mathematicians in order to explain the growth of
quasicrystals in terms of tilings. In this model, the atoms are
located at the vertices of the tiles.

Both models can be obtained from the following construction. Let a
``superlattice'' $\Gamma$ of full rank be given
in the superspace $\real^n$, and the
superspace be presented as the direct sum
$\real^n=\real^k\bigoplus\real^{n-k}$,
where $\real^k$ is the physical subspace.
Let
$p:\real^n\rightarrow \real^k$ and $q:\real^n\rightarrow
\real^{n-k}$ be the natural projections.
Fix a
finite $(n-k)$-polyhedron $D\subset
\real^{n-k}$ and consider the ``tubular $D$-neighborhood''
$D_q=q^{-1}(D)\subset \real^n$ of the physical subspace
$\real^k\subset\real^n$.
Assume that the boundary of the polyhedron $D$ is disjoint from
$q(\Gamma)$, or equivalently, $\partial D_q\cap\Gamma=\varnothing$.
Then the set of points
$$p(\Gamma\cap D_q)\subset\real^k$$
is quasiperiodic in the sense of the definition given above.

By taking a certain polytope
decomposition of the space $\real^n$ associated with the lattice
$\Gamma$ and the polyhedron $D$, one obtains a quasiperiodic
tiling of $\real^k$
whose tiles are the intersections of $\real^k$ with
the $n$-cells of the decomposition (see survey
article \cite{Pi}).

Very interesting examples of nontrivial
symmetry groups come from the superspace $\real^4$ endowed with
the Minkovski metric and
a lattice $\Gamma\cong\integer^4$ such that the
physical subspace $\real^2$, which is spacelike
(\emph{i.e.}\ Euclidean), is
invariant under some lattice-preserving mapping from the group
$O(3,1)$.

The superspaces $\real^{l,m}$, where $l+m=n$, might also appear in
interesting cases.

\subsection{Quasiperiodic functions in the theory of conductivity}

Here we describe the situation that is the main motivation for
our topological and dynamical theory.

For every single
crystal normal metal, we have a lattice $\Gamma$ in the physical
space $\real^3$. However, our geometrical constructions will live in
a completely different space, namely the $3$-torus of
quasimomenta $\torus^3$, which
is the quotient space of the \emph{dual} $3$-space $(\real^3)^*\cong\real^3$
by the dual
(reciprocal) lattice $\Gamma^*\cong\integer^3$. The ``Bloch'' states of quantum
electrons are parameterized by pairs $(m,p)$,
where $p$ is a point in the space of
quasimomenta, $p\in \torus^3=\real^3/\Gamma^*$, and $m$ is a
natural number, which is the index of
a branch of the dispersion relation
$f(p)=\epsilon_m(p):\torus^3\rightarrow \real$. In what follows
we will always deal with just one branch only, so we
drop the index $m$ in the notation. We assume that $f(p)=\epsilon(p)$ is a Morse
function on the $3$-torus or, in other words,
a three-periodic Morse function on the
covering Euclidean space $\real^3$. At zero temperature all
electrons occupy the ``Dirac sea'' $\epsilon(p)\leqslant \epsilon_F$
where the ``Fermi energy'' $\epsilon_F$ depends on the number of free
electrons in the metal. We assume that $\epsilon_F$ is a regular value
for the Morse function $f=\epsilon(p)$. At low temperatures we
are dealing only with
``excited'' electrons nearby the
Fermi level $\epsilon(p)=\epsilon_F$.

The Fermi level looks geometrically as a
two-dimensional surface $M_F\subset \torus^3$
in the space of quasimomenta. This
surface is nonsingular and homologous to zero in the $3$-torus. Let
us assume that it is connected.

The \emph{topological rank} $r(M_F)$
of the Fermi surface is defined as the rank of the image
of the first homology group of $M_F$ under the mapping
$i_*:H_1(M_F,\integer)\rightarrow H_1(\torus^3,\integer)\cong\integer^3$
induced by the inclusion $i:M_F\hookrightarrow\torus^3$.
Since $i_*(H_1(M_F))$ is a sublattice in $\integer^3$,
we always have $r\in\{0,1,2,3\}$.

For
example, the topological rank of the Fermi surface of
lithium is equal to zero (in this case, the Fermi
surface looks like a topological $2$-sphere), whereas
it is equal to three for
copper, gold, platinum, and some other noble metals. For gold, for
example, the genus of the Fermi surface is equal to four.

The problem that we will consider is most difficult when the topological
rank of the Fermi surface is maximal possible, \emph{i.e.}, equal to three.
One can easily show that the genus of the Fermi surface must be
greater than or equal to the topological rank.

Interesting dynamical phenomena occur in the presence of a
magnetic field. In the semiclassical approximation, an electron,
which is considered as a point in the space of quasimomenta, moves
along constant energy lines in the plane $\real^2_{B,p_0}$
orthogonal to the magnetic field $B$ and passing through the
initial position $p_0$ of the quasimomentum.

One may say that this is a
Hamiltonian system on the $3$-torus of quasimomenta with Poisson
bracket $$\{p_j,p_l\}=\frac ec\,B_{jl}=\frac ec\,\varepsilon_{jlq}B^q$$ and
Hamiltonian $f=\epsilon(p)$:
$$\frac{dp_j}{dt}=\{p_j,\epsilon(p)\},$$
so
the motion preserves the energy and a linear Casimir of the Poisson
bracket. The level sets of this Casimir are planes orthogonal to the magnetic
field. The trajectories can be treated as  leaves of the
Hamiltonian foliation on the Fermi surface given by the
equation $\omega=0$ where $\omega$ is the following closed $1$-form:
$\omega=\sum_jB^jdp_j|_{M_F}$.

According
to the ``strong magnetic field limit'' principle worked out by
I.\,Lifshitz, M.\,Azbel, M.\,Kaganov, and V.\,Peschanski in early 1960s,
all essential properties of the electrical conductivity in the presence of
a reasonably strong uniform magnetic field $B$ should follow from the
structure of the
dynamical system on the Fermi surface described above (see
\cite{LAK,Ab,Ki,Zi}). For ordinary normal metals (like gold, for
example) one may use this approximation for magnetic fields strong
enough in the human sense (like
$1\,\mathrm{Tesla}<|B|<10^3\,\mathrm{Tesla}$ for low
temperatures; recall: $1\,\mathrm{Tesla}=
10^4\,\mathrm{Gauss}$). If the magnetic field is too strong, then
the semiclassical approximation will not be valid.
If the magnetic field is too weak, then the electron quasimomentum drift will
be too slow, and the distance that the quasimomentum passes for the
characteristic time of the electron free motion will become
insufficient to affect the observable conductivity.

However, in 1960s the study of the just mentioned dynamical system
was only started. Some conceptual mistake
was then made in \cite{LP} and further investigation
was stopped, and resumed only many years later in works
\cite{N,Z,N1,D,D1,NM,NM1,D2,D3,DM,MN,MN1,MN2}.

What is crucial for us here is following:

\medskip
\centerline{\parbox{0.8\textwidth}{%
the electron trajectories coincide with connected components of
the level curves $\epsilon(p)=\epsilon_F$ of the function
$\epsilon$ restricted to the planes orthogonal
to the magnetic field $B$; in other words, they are connected
components of the level curves of functions that form a family of
related quasiperiodic functions with three quasiperiods. }}
\medskip

In work \cite{N2} an extension of these studies to
a larger number of quasiperiods was started. In particular, some
new  ideas and results   were formulated for the case $n=4$. The
present  work contains the first complete proof of those (properly corrected)
statements. The proof is based on the topological technique developed in
\cite{D2,D3}.

Let a constant Poisson bracket $B_{jk}$ of rank two be given
on the $n$-torus. Then every Hamiltonian
$f(p)=\epsilon(p):\torus^n\rightarrow \real$ defines a Hamiltonian
system whose trajectories are exactly the level lines
$\epsilon(p)=\const$ of the restriction
of the function $\epsilon$ to the planes $\real^2_{B,a}$
defined as follows.
There exist exactly $n-2$ independent linear
Casimirs $K_1,\dots,K_{n-2}$, $K_j(p)=K_j^lp_l$, such that
$\{p_s,K_j\}=0=K_j^lB_{sl}$. We put
$$\real^2_{B,a}=\{K_1=a_1,\dots,K_{n-2}=a_{n-2}\},$$
where $a=(a_1,\dots,a_{n-2})$.
So our trajectories
are exactly the levels of quasiperiodic
functions on the two-planes $\real^2_{B,a}$, which form
the family of descendants of the $n$-periodic function
$\epsilon(p)$. They depend on the constants $a_1,\dots,a_{n-2}$.

Topological study of this problem is the central part of this
article.

Modern experimental technology allows to construct surfaces with
a variety of prescribed small fluctuations. In particular, it is
possible to make  a quasiperiodic construction with any number of
quasiperiods. It presents us a two-dimensional weak quasiperiodic
electric potential $V(x,y)$. In a strong magnetic
field $B$ electrons move along the surface. After averaging we
obtain a slow motion along the level curves $V(x,y)=\const$. These
studies, experimental and theoretical, were done originally for
periodic potentials with $n\leqslant 2$ periods only (see
\cite{B}), but it was pointed out in work \cite{M} that
quasiperiodic potentials can also appear here; new predictions
were made for quasiperiodic cases with three and four quasiperiods
based on the topological results obtained in a series of works
of the present authors (\cite{N2,D3,MN2}).

\section{Topology and dynamics of quasiperiodic functions on the
plane: the case of three quasiperiods. The electrical conductivity in
metals}

We address the following general question. How may the level lines
$\varphi=\const$ of a quasiperiodic function $\varphi$ on the plane
with $n$
quasiperiods look like? In a generic situation, such a level line
is a one-dimensional submanifold of $\real^2$, i.e.\ a union of
curves. We will call these curves ``trajectories'' because in our studies they
have been appearing as semiclassical electron trajectories
on the Fermi surface in the presence of a magnetic field since early 1980s when
this problem was posed in work \cite{N} as a problem of
topology and dynamical systems. It corresponds to the case of three
quasiperiods only. Some of those curves may be closed in $\real^2$
(compact) and others nonclosed in $\real^2$ (open). Let us ask the
following questions.

\begin{description}
\item[Question 1] Is the size of the compact trajectories uniformly bounded (for a
fixed level of $\varphi$)?
\item[Question 2] Do the open trajectories have some nice asymptotic behavior?
\end{description}

The first results were
obtained in work \cite{Z}. It became clear in the second half
of 1980s  that the proper form of Question~2 is the following:
does any open trajectory has a ``strong asymptotic direction'' in
$\real^2$, \emph{i.e.}, lie in a strip of uniformly bounded width
and passes through the strip
``from $-\infty$ to $+\infty$''? This specification of the problem
was made in article \cite{N1}. In work \cite{D} the
results of \cite{Z} were improved accordingly to the new formulation of
the problem. The important breakthrough was made in work
\cite{D1}, but for a long time there was no applications.
Physical applications were found later in works~\cite{NM,NM1}.

In the physically important  case $n=3$, the positive answer to
our Question~1 follows easily from a quite elementary argument.
For $n>3$ it is more difficult, and it will be discussed later.
Question~2 is highly nontrivial already for $n=3$.
As mentioned above, the asymptotic behavior of open electron
trajectories, \emph{i.e.}, of open connected components
of a level line of a quasiperiodic
function with three quasiperiods, was studied in
\cite{Z,D,D1}. It became finally
clear after work \cite{D1} that for the family of
related quasiperiodic function corresponding to
a ``typical'' direction of the
magnetic field (which is regarded as the direction of a
plane $\real^2\subset \real^3$),
either their level lines do not have open components at all or
the open components all have a strong asymptotic direction.
 The
latter means that each open curve has a parametrization
$\gamma(t)$ such that the following holds for some nonzero two-vector
$(x_1,y_1)$:
\begin{equation}\label{strong}
\gamma(t)=(x_0,y_0)+t\cdot(x_1,y_1)+O(1).\end{equation}
Below we will explain precisely
what  `a typical direction' means,
and provide references to the papers containing proofs of
the corresponding results. In the first work \cite{Z}
completed in the note
\cite{D}, this type of result was obtained for the special case of
small perturbations of a magnetic field having ``rational''
direction. We shall return to this special case in the next section
where we discuss the
quasiperiodic functions with four quasiperiods.

Applications of this studies to explaining the electrical
conductivity in a strong magnetic field are presented in
works \cite{NM,NM1}. They are based on the results of the Lifshitz
school of 1960s. Physicists calculated the contribution of
individual trajectories of simple types to the conductivity
tensor. These calculations have become a part of textbooks (see
\cite{LAK,Ab}). In the case when all trajectories are
compact the conductivity
components orthogonal to the magnetic field $B$ decrease as $|B|^{-1}$ or
$|B|^{-2}$ when $|B|$ grows while the direction of $B$
remains fixed. Some special examples of open trajectories lying in
finitely wide strips were found at the same time and their
contribution to the conductivity was
calculated. As pointed out in \cite{NM,NM1},
one can easily extend the just mentioned
calculation to the case of general trajectories of the same type.
The projection to the plane orthogonal to $B$
of the part of the conductivity tensor contributed by such a trajectory
has two eigenvalues one of which is zero and the other nonzero.
Since the contribution of closed trajectories tends to zero
when $|B|$ grows, the observable conductivity tensor
for a strong enough $B$ will depend on open trajectories only.

\smallskip
\centerline{\parbox{0.8\textwidth}
{However, the observable physical
conductivity tensor is formed by the contributions of all electron
trajectories as the sum of them. What conclusion about this tensor
can be made from the qualitative dynamical properties of
that system, which was defined on the quantum level?}}
\smallskip

In order to obtain a nontrivial new physical result, one needs
more than the theorems explicitly formulated in~\cite{D1}.
But luckily an additional crucial property also holds
for our dynamical system,
and this can be extracted from the proofs of the main theorems
of works \cite{Z,D1}. The property, which we call
\emph{topological resonance}, implies the following for the behavior
of trajectories.

\smallskip
\centerline{\parbox{0.8\textwidth}{For a ``typical'' family $\{\varphi_a\}$
of related quasiperiodic functions $\varphi_a(p)=\epsilon(p)|_{\real^2_{B,a}}$,
the strong asymptotical direction $\eta_B$ of noncompact trajectories is the
same for all trajectories.
Moreover, there exist an integral two-plane $\mu\subset
\real^3$ (\emph{i.e.}, a plane generated by two reciprocal lattice vectors,
$\mu\cap\integer^3\cong\integer^2$)
such that $\eta_B$ has the direction of the intersection of $\mu$
with the plane orthogonal to the magnetic field:
$$\eta_B\in\mu\cap \real^2_B.$$
This integer plane $\mu$ is
locally rigid, \emph{i.e.}, it remains unchanged under small
variations of the direction of the magnetic field.}}
\smallskip

The topological resonance property of our dynamical system
makes possible serious applications. It was missed in the classical works
of physicists, and a conceptual mistake was made in \cite{LP},
where calculations led to a result contradicting to this
property. This mistake was revealed and corrected only in works
\cite{NM,NM1,MN}.

For a strong enough $B$, the direction $\eta_B$ is a
zero eigenvector of the projection of the conductivity tensor to the
plane orthogonal to the magnetic field. The
integral plane $\Pi\subset \real^3$ is
directly observable by measuring the zero eigenvector $\eta_B$ for two
or more magnetic fields $B$ close to each other.

We refer the reader to recent article \cite{MN} for a more detailed
physical discussions.

Let us describe the picture topologically. Consider all our objects
in the universal covering space $\real^3$ with the reciprocal lattice
$\Gamma^*=\integer^3\subset \real^3$ and the three-periodic Fermi surface
$$\widehat{M_F}=\nu^{-1}(M_F)\subset \real^3$$
covering the compact one
$$M_F\subset
\torus^3.$$
(Recall that $\nu$ stays for the standard projection $\real^3\rightarrow\torus^3$.)
The electron trajectories in the covering space are
connected components of the intersections of the
three-periodic Fermi surface with planes orthogonal to $B$. Let
$M_0(B)$ be the closure of the union of all compact trajectories,
and let $L(B)$ be the closure of its complement in the Fermi surface:
$$L(B)=\overline{M_F\setminus M_0(B)}.$$

Let $L_l(B)$ be the connected components of $L(B)$.
In the typical case, $L(B)$ is a compact two-manifold with boundary
$$\partial
L(B)=\bigcup_{l,s} \beta_{ls},$$
where
$$\partial L_l(B)=\bigcup_s\beta_{sl}.$$
All boundary  curves $\beta_{sl}$ are
saddle connection cycles. In the typical case, we may assume that
every cycle $\beta_{sl}$ joins a saddle critical point to itself, since
all the other cases have
positive codimension in the appropriate functional space.
(In particular, we assume that there
is no rational linear dependence between the components of $B$,
and that the Hamiltonian foliation defined by $\omega=\sum
B^idp_i|_{M_F}=0$ has only Morse singularities and
does not have saddle connections between different saddles.)

The part $M_0$ of the Fermi surface
can be presented as the union of ``cylinders'' $Z_q$,
$M_0=\bigcup_q Z_q$,
whose interior consists of regular compact trajectories and
``bases'' are either saddle connection cycles or
isolated points (centers). There are finitely many such
cylinders, and they are obviously compact.
This immediately implies a positive
answer to Question~1 posed in the beginning of this section:

\begin{center}the size of all compact trajectories is uniformly bounded.
\end{center}

We call the
pieces $L_l$ of the Fermi surface the \emph{carriers of open
trajectories}. By construction, every open trajectory
(in $\torus^3$) is contained in one of the carriers,
and, in the generic case, is everywhere dense in it.
Let $D^2_{ls}\subset \real^2_{B,a}$ be planar two-discs
orthogonal to the magnetic field such that $\partial
D^2_{ls}=\beta_{ls}$. We define the ``closure'' $N_l$ of every carrier
$L_l$ as follows: $$N_l=M_l\cup\Big(\bigcup_{s} D^2_{ls}\Big).$$
By construction we also have
$$N_l\bigcap N_{l'}=\varnothing$$
for $l\neq l'$.

We call
our system {\it stable topologically completely integrable} if
the genus of each surface $N_l $ is equal to one, and this
picture is stable under arbitrary small enough perturbations of the magnetic field.

We call the system {\it chaotic} if the genus of some $N_l$ is greater
than one. According to the main theorem of \cite{D1}, the latter
situation is always topologically unstable.

According to \cite{NM,NM1,MN}, what is important for physical
applications, is the following \emph{topological resonance}
property of our system. In the stable topologically completely
integrable case, all the closures $N_l$ of the carriers of open
trajectories have the same up to sign nonzero homology class:
$$[N_l]=\pm \mu\in H_2(\torus^3,\integer),\quad\mu\neq 0,$$
which is an indivisible element of the group
$H_2(\torus^3,\integer)\cong\integer^3$. The number of the tori
$N_l\subset \torus^3$ is even because the sum of their homology classes is
equal to the class of Fermi surface, which is zero.
(Note that every homologically nontrivial connected closed nonselfintersecting
two-manifold $M\subset \torus^3$ always represent an indivisible
homology class. Any two such submanifolds with empty
intersection represent the same homology class up to
sign.) The class $\mu\in H_2(\torus^3,\integer)\cong\integer^3$
is presented by three relatively prime integers $\mu(B)=(m_1,m_2,m_3)$.

The integral vector $\mu(B)$ remains
unchanged under small perturbations of $B$. Therefore, there
is an open set on the sphere $S^2$ with the same $\mu(B)$.
This set as well as the integral vector $\mu$ is
an observable characteristics of our system,
and it can be found experimentally by measuring the
conductivity tensor in the presence of strong
enough magnetic fields having generic directions.
The stable topologically integrable
case occurs for all directions $B/|B|\in S^2$
of the magnetic field from an everywhere dense open subset of
$S^2$. It was proved
in \cite{D2,D3} (by two different methods)
that this picture may be not valid
for directions $B/|B|\in S^2$ of the magnetic field from
a nonempty subset whose codimension is at least one.

In terms of quasiperiodic functions, we can say that ``non-typical''
functions $\varphi_a(p)=\epsilon(p)|_{\real^2_{B,a}}$ with three quasiperiods,
\emph{i.e.}, such that open connected
components of their level sets don't have a strong asymptotic direction,
all lie in a subset that has codimension one (in some
natural sense). Examples of level lines with chaotic behavior in
the case $n=3$ were constructed in \cite{D2}. We call such level lines
\emph{strongly chaotic trajectories}. Interesting attempts were
made in order to find physical properties of the conductivity in these
cases. For some special examples it was done in work \cite{DM}
but in general the stochastic properties of these trajectories are
unknown. A.\,Maltsev formulated the following conjecture.

\begin{conjecture}
The contribution of strongly chaotic trajectories to the
conductivity tensor tends to zero when $|B|$ grows (remaining in
a reasonable range), which includes the conductivity in
the direction of the  magnetic field itself.
\end{conjecture}

Previously, S.\,Tsarev constructed a ``weakly chaotic'' example
(unpublished, see work \cite{D2}). In his case, there is a rational dependence
between the components of $B$, and there is just one carrier of
open trajectories in our sense, which coincides with the Fermi surface.
However, the closure of any trajectory in $\torus^3$ is not the whole
surface, but just a half of it, which is homeomorphic
to a $2$-torus with two holes. The holes are not homologically
trivial in $\torus^3$, so they are not regarded as closed in $\real^3$.
In Tsarev's example, the nonclosed level lines still have an
asymptotic direction in a weaker sense,
$$\gamma(t)=(x_0,y_0)+t\cdot(x_1,y_1)+o(t),$$
but the projection of any trajectory to a straight line perpendicular to $(x_1,y_1)$
is unbounded.

It is interesting to look at the behavior of trajectories in the
special (nongeneric) case of the Fermi surface
$$\epsilon(p)=\cos(p_1)+\cos(p_2)+\cos(p_3)=0.$$
Examples of this type were investigated numerically and
analytically in works \cite{D2,RdL}. There are chaotic trajectories for the set of
magnetic fields whose Hausdorf dimension is presumably equal to some
$\alpha$ with
$$1<\alpha<2.$$
There are many (in fact, infinitely many)
different stable topologically completely
integrable zones on the sphere $S^2$ having different integral
characteristics $\mu(B)\in\integer^3$. We
call this type of examples \emph{generic symmetric
levels}, see below.

We conjecture the following.

\begin{conjecture}
(i) For a generic connected smooth two-manifold $M_F\subset
\torus^3$ homologous to zero, the set of chaotic directions of
the magnetic field has Hausdorf dimension less than one in $S^2$.

(ii) For a generic $1$-parametric smooth family $M_{F,t}\subset \torus^3$
of such Fermi surfaces this set has Hausdorf
dimension less than two.
\end{conjecture}

A detailed investigation of this problem containing the proofs of
all topological statements needed for physical
applications found in \cite{NM,NM1} is performed in
\cite{D2,D3}. Special attention is paid there to one-parametric families of
Fermi surfaces that are levels of the same Morse function
$f:\torus^3\rightarrow \real$:
$$M_c=\{f(p)=c\}.$$
It is proved that, for any $B$ from a stability zone,
open trajectories live on the levels $M_c$ from a connected interval $c_1(B)\leqslant
c\leqslant c_2(B)$ of the real line. For a $B$ away from the stability
zones, the strongly chaotic
behavior might appear only on a single level $c(B)\in\real$.
As a corollary we obtain the following result:

\smallskip
\centerline{\parbox{0.8\textwidth}{For a function
f=$\epsilon(p)$ with symmetry $\epsilon(p+p_0)=-\epsilon(p)$,
where $p_0\in\torus^3$ is some shift,
strongly chaotic trajectories cannot appear on the levels $M_c$ with
$c\neq0$ because otherwise they must appear on $M_{-c}$, too,
for the same $B$, which is impossible.
In such a case we call the level $c=0$ a \emph{generic symmetric level}.}}
\smallskip

According to
our conjecture, the Hausdorf dimension of the set of $B/|B|\in S^2$ for which
the strongly chaotic behavior occurs on such a level is equal to some $\alpha<2$.

The surface $\sum_{j=1}^{j=3} \cos(p_j)=0$ gives an example of a generic
symmetric level with $p_0=(\pi,\pi,\pi)$.

Some more details about chaotic
trajectories and stability zones for the case of three
quasiperiods will be given below. They will be needed for
proving our main
result about quasiperiodic functions with four quasiperiods (see the next
section).

\section{The stable topological complete integrability for $n=4$
quasiperiods}

Let us consider now the case of $n=4$ (or more) quasiperiods.
For every direction $\Pi$ of two-planes in $\real^n$ (\emph{i.e.}, a two-dimensional
vector subspace $\Pi\subset\real^n$), the
original $n$-periodic Morse function
$f:\torus^n\rightarrow \real$
defines a family of descendants $\{\varphi_a(y)\}$ on the family
affine two-planes $\real^2_{\Pi,a}\subset
\real^n$ having direction $\Pi$.

We call the level $\{f=c\}$
of the function $f$ \emph{topologically completely
integrable} (\emph{TCI}) for the direction $\Pi$ if, for each
$\varphi_a$ from the family, all regular connected components of the level
line $\varphi_a(y)=c$ are either compact or have a strong asymptotic
direction. We call this level \emph{stable
TCI} if this property remains unchanged under small
perturbations of the function $f$ and the direction $\Pi$,
which is a point in the Grassmanian manifold $G_{n,2}$.

We say that
the Stable TCI level satisfies the {\it topological
resonance} condition (for given $\Pi$) if there exists an integral hyperplane
$\mu\subset \real^n$, $\mu\cap\integer^n\cong\integer^{n-1}$,
such that all open regular trajectories have
the same asymptotical direction $\eta_\Pi$ that coincides with
the direction of the straight line $\mu\cap\Pi\cong\real$.
Since $\mu$ is integral, it must remain unchanged under small perturbations of
anything.

Let us point out that even the ``trivial case'' $n=2$ is
meaningful (as a subject of the elementary Morse theory on the
$2$-torus): for a generic double-periodic function on the plane
there exists a level $f=c$ with a connected component presenting a nontrivial
indivisible
homology class $\mu\in H_1(\torus^2,\integer )$.
All other components of every level are
either homologically trivial or homologous to $\pm\mu$. For
Morse functions with exactly four critical points and critical values
$c_0<c_1<c_2<c_3$, all levels $f=c$, with $c_1<c<c_2$, have exactly two
connected components, whose homology classes are $\pm\mu$. All other nonsingular
levels are either compact or empty.

\begin{description}
\item[Question]
Consider famous real nonsingular quasiperiodic
finite-gap solutions of the KdV equation
$$u(x,t)=2\partial_x^2\log\Theta (xU+tV+U_0) + c_{\Gamma}$$ with
an arbitrary number of quasiperiods (or gaps). Are the levels
$u(x,t)=c$ always stable topologically completely integrable or
they can be chaotic? How to find their strong asymptotic
direction and their integer-valued characteristic $\mu$?
\end{description}

P.\,Grinevich pointed out to us that, for real smooth finite-gap
solutions of the KdV equation,
the $n$-periodic function $f=2\partial_U^2 \log
\Theta(\eta_1,\dots,\eta_n)+c_{\Gamma}$ on the $\eta$-space
is always a Morse
function on the real $n$-torus with $2^n$ critical points, simply
because it can be reduced to the form
$f=\sum_{j=1}^{j=n}\alpha_j\sin x_j$ by a diffeomorphism of the
torus isotopic to the identity. There is a canonical  lattice in the
$\eta$-space generated by the so-called $a$-cycles which are the
real finite gaps of the 1D Schr\"odinger operator on a hyperelliptic
``spectral'' Riemann surface $\Gamma$. The real constants
$\alpha_j$ depend on the spectrum. As a conclusion, we get the
following:

\smallskip
\centerline{\parbox{0.8\textwidth}{For
generic real nonsingular two-gap solutions of the KdV equation, there
exist a critical value $c_{\mathrm{cr}}$ such that all constant speed
levels $u(x,t)=c$,
$$c_1=c_{\Gamma}-c_{\mathrm{cr}}<c<c_{\Gamma}+c_{\mathrm{cr}}=c_2$$
are periodic
perturbations  of a family of straight lines with integral
direction $m_1:m_2$ on the plane
with lattice. This direction is locally rigid, but globally depends
on the constants $\alpha_j$. All other levels are either compact
or empty. We call it the \emph{topological speed} of the solution.}}
\smallskip

The computational studies of this problem for finite-gap solutions
are now being investigated numerically.

We concentrate now on the case $n=4$ quasiperiods. Work
\cite{N2} presents an idea of the proof that, for every generic Morse
function $f$ and a noncritical generic level $f=c$, there exists an open
everywhere dense set of two-plane directions
$\Pi\in G_{4,2}$ for which the level $f=c$
is stable TCI. The proof of this theorem
requires the use of some extension of the results of work
\cite{D3}. Here we make some corrections to the statement
of~\cite{N2} and provide a complete proof.

We believe that a generic level is stable TCI for all directions $\Pi$
from a subset $S\subset G_{4,2}$ whose measure is full in $G_{4,2}$. However,
we don't have an idea how to prove this conjecture. For $n>4$
nothing like that is expected.

Now we start a detailed investigation of the case $n=4$. Even
Question~1 of the previous section presents a difficulty here.
It is possible that a
single level set of a quasiperiodic function with four
quasiperiods contains a family of compact components without an
upper bound for their size. An example can be constructed easily.

However, we are going to show that there is an open everywhere
dense open set of
quasiperiodic functions $\varphi$ with four quasiperiods such that
the level lines $\varphi=\const$ have the same qualitative
behavior as those in the typical case of three quasiperiods.
The precise formulation is as follows.

\begin{theorem}\label{th1}
There exists an open everywhere dense
subset $S\subset C^\infty(\torus^4)$ of $4$-periodic
functions $f$ and an open everywhere dense subset
$X_f\subset G_{4,2}$ depending on $f$ such that
any level $M^3_c=\{f=c\}$ of $f$ is stable TCI (or does not
contain open trajectories at all) for any $\Pi\in X_f$.

Moreover, for any regular open trajectory, the remainder
term $O(1)$ in (\ref{strong}) as well as
the diameter of any compact trajectory are bounded from above by a constant $C$
not depending on the affine plane $\real^2_{\Pi,a}$ containing the trajectory,
provided that the level $c$ and the direction $\Pi\in X_f$ are fixed.
\end{theorem}

Let us make a remark about notation and terminology.
Once we switched to the case of four quasiperiods,
our problem is no longer relevant to the discussed above physical model
of conductivity in normal metals in the presence of a magnetic field.
So we change the notation for the coordinates in $\real^n$ from $p_l$, which
we used for quasimomentum, to more customary, $x_l$, $l=0,1,2,3$,
and don't longer think
of the lattice $\integer^4\subset\real^4$ as the one dual to
some physical lattice. We think of the ``magnetic field''
$B$ as a linear mapping from $\real^4$ to $\real^2$ (or from $\real^3$
to $\real$ in the $n=3$ case) such that $\Pi=\ker(B)$.
Thus, by $\real^2_{\Pi,a}$ we mean the two-plane $B^{-1}(a)$,
where $a\in\real^2$. In the case $n=3$ we may also think of $B$ as a vector
perpendicular to the plane $\Pi$.
However, we keep calling connected components of
the intersections of $M_c^3$ with the two-planes $\real^2_{\Pi,a}$
trajectories, just for briefness.

We start by recalling results of~\cite{D2,D3} for the
three-dimensional case in the form needed to prove our theorem. Let
$B:\real^3\rightarrow\real$ be a linear function of irrationality
degree three, \emph{i.e.}, of the form $B(x)=B_1x_1+B_2x_2+B_3x_3$,
where $B_1$, $B_2$, $B_3$ are reals linearly independent over $\integer$, and
let $f:\torus^3\rightarrow\real$ be a generic smooth function. By
`generic' we mean that $f$ does not satisfy certain conditions
that have codimension $\geqslant1$. However, we shall pay
attention to codimension one singularities as, in order to deal
with the four-dimensional case, we are going to consider
one-parametric families of three-dimensional pictures.

We abuse notation by using the same letter $f$ for the lift of $f$
to the covering $3$-space $\real^3$. We use notation $M^2_c$
for the level set $f^{-1}(c)$ in
$\torus^3$ and $\widehat{M^2_c}$ for its cover in $\real^3$. By
$\gamma_{a,c}$ we denote the whole intersection of $\widehat{M^2_c}$
with the plane $\real^2_{\Pi,a}=B^{-1}(a)$. So, the trajectories that we are studying are
regular connected components of $\gamma_{a,c}$ or their projections to $\torus^3$.

First of all, consider closed trajectories on $M^2_c$. Notice
that, since we assumed $B$ to be of maximal irrationality degree,
a trajectory in $\real^3$ is closed if and only if so is its image
in $\torus^3$. Without the assumption on $B$ this may be not true,
since a closed trajectory in $\torus^3$ may then be
non-homologous to zero, in which case its cover in $\real^3$
consists of infinite ``periodic'' trajectories treated as ``open''
in the physical applications.

For a generic $f$, compact trajectories on every $M^2_c$ form
finitely many cylinders whose bases are either saddle connections
or extrema of the restriction $B|_{M^2_c}$. Obviously, the length
of compact trajectories is bounded from above by some constant.

Let $U$ be the set of $c$ such that $\gamma_{a,c}$ has unbounded
connected components for some $a$. In other words, $c\in U$ if and
only if $M^2_c$ contains open trajectories that are not saddle
connections. The following picture, which was sketched in
the previous section, can be extracted from work
\cite{D3}:

\smallskip
\noindent\hskip 0.1\textwidth{\parbox{0.8\textwidth}{The
set $U$ is either a closed interval, $U=[c_-,c_+]$, or
just one point, $U=\{c_0\}$.

If $U=[c_-,c_+]$ is a nontrivial interval, then for any $c\in U$,
there is a (unique) family of two-tori $\torus^2_{c,1},\ldots,
\torus^2_{c,2k}$ (with $k$ depending on $c$)
imbedded in $\torus^3$ such that
\begin{enumerate}
\item Each $\torus^2_{c,i}$ consists of the closure of some
open trajectory on $M^2_c$ and a few (may be zero) planar
disks perpendicular to $B$;
\item Every open trajectory is contained by whole in one of
the tori $\torus^2_{c,i}$;
\item All the tori $\torus^2_{c,i}$ define the same up to sign
nonzero homology class $\mu$ in $H_2(\torus^3,\integer)$;
\item For all but finitely many $c$, the tori $\torus^2_{c,i}$
are pairwise disjoint, and, in this case, a sufficiently
small variation of $c$ causes small deformation of the tori.
For the exceptional $c$'s they can be
made disjoint by a small perturbation. At such a $c$, a couple of tori
is born or killed.
\end{enumerate}
All the picture is stable in this case, which means that after a
small enough perturbation of $f$, the interval $U$ and the family
of tori $\torus^2_{c,i}$ are perturbed slightly. In particular,
the homology class $\mu$ fixed.}}\hfill\\

\begin{remark}
In the setting of papers~\cite{D2,D3},
the function $f$ was assumed to be fixed
and the point of concern was the dependence of the behavior of
our dynamical system on the magnetic field $B$ and on the
level of the function $f$. The stability
of the whole picture under small perturbations of the function $f$
was not discussed. However, the arguments of those works can be
easily modified in order to prove such stability. Indeed,
one of the key observations in~\cite{D2,D3} is that, locally,
the qualitative behavior of the trajectories (including the existence
of a strong asymptotic direction) depends only
on finitely many parameters, which are certain critical values of
the ``height'' function $B(x)$ restricted
to the surface $\{f=\const\}$. (For instance,
the existence of strongly chaotic examples
was proved in~\cite{D2} by specifying the combinatorial
structure of the surface and particular values of the
parameters.) It is easy to see that those parameters
behave nicely under small perturbations of $f$,
so extending the arguments of~\cite{D2,D3}
to this, more general type of perturbations requires
almost no additional work.
\end{remark}

\begin{figure}[ht]
\centerline{\epsfig{file=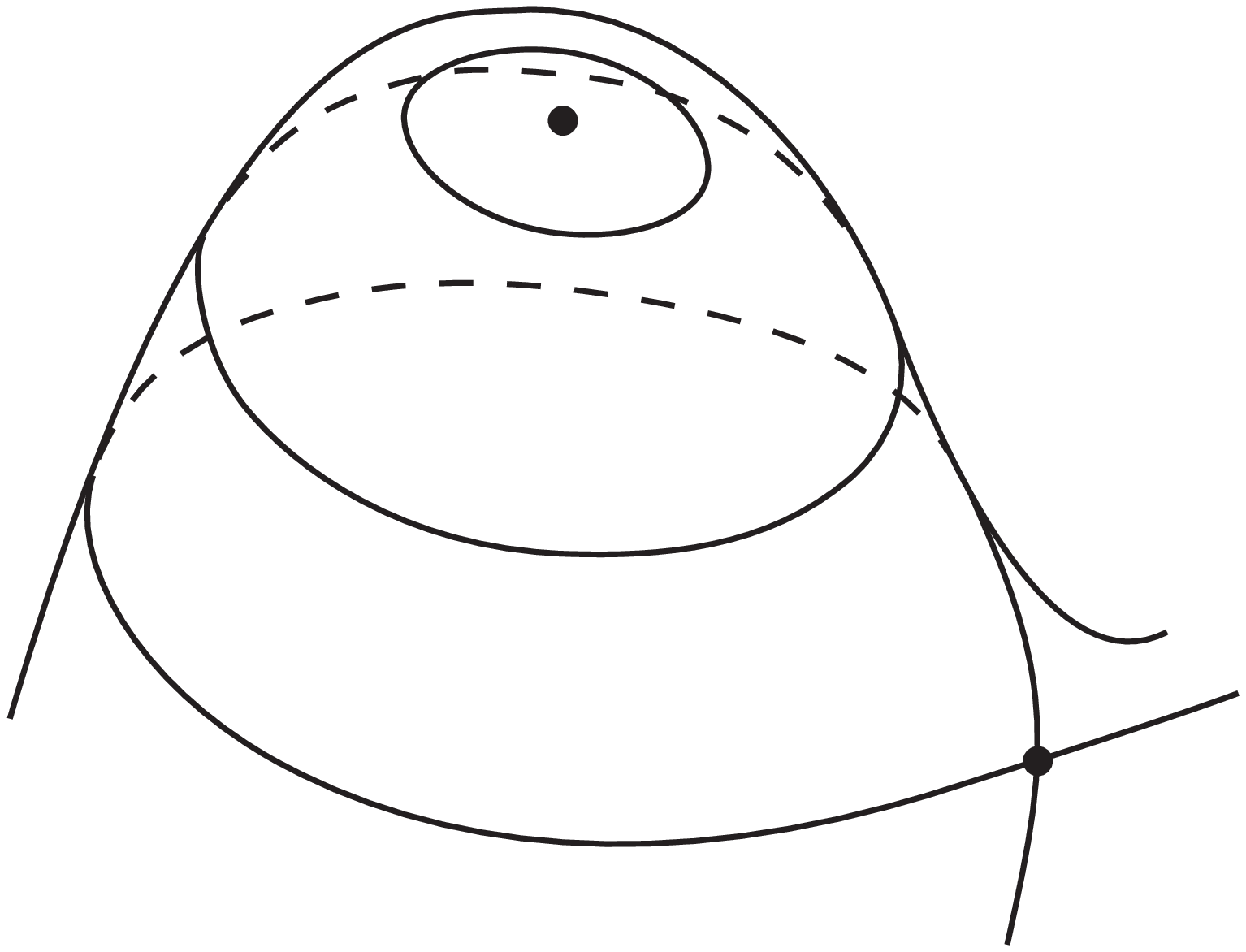,height=100pt}\hskip1cm
\epsfig{file=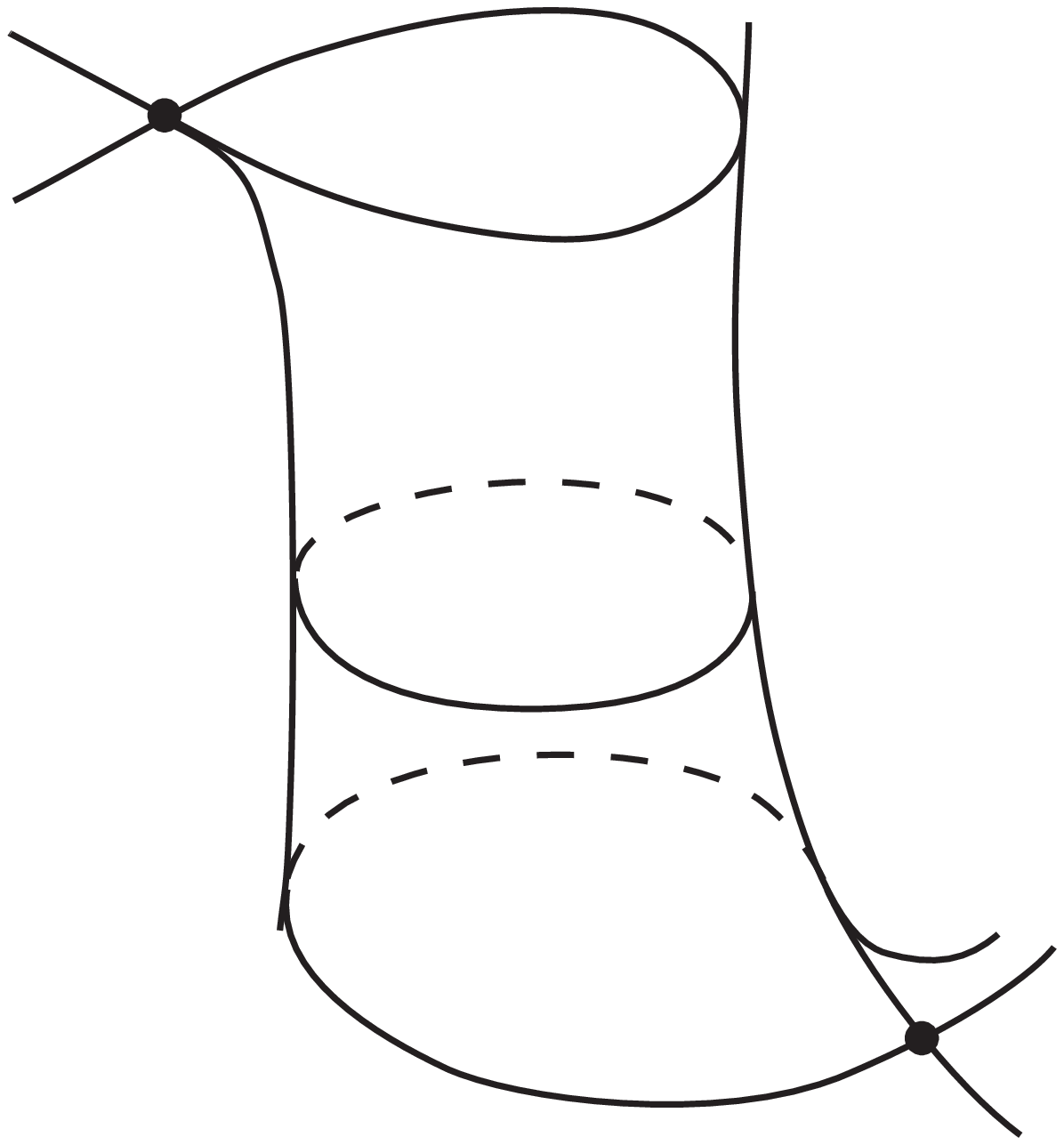,height=100pt}\hskip1cm\epsfig{file=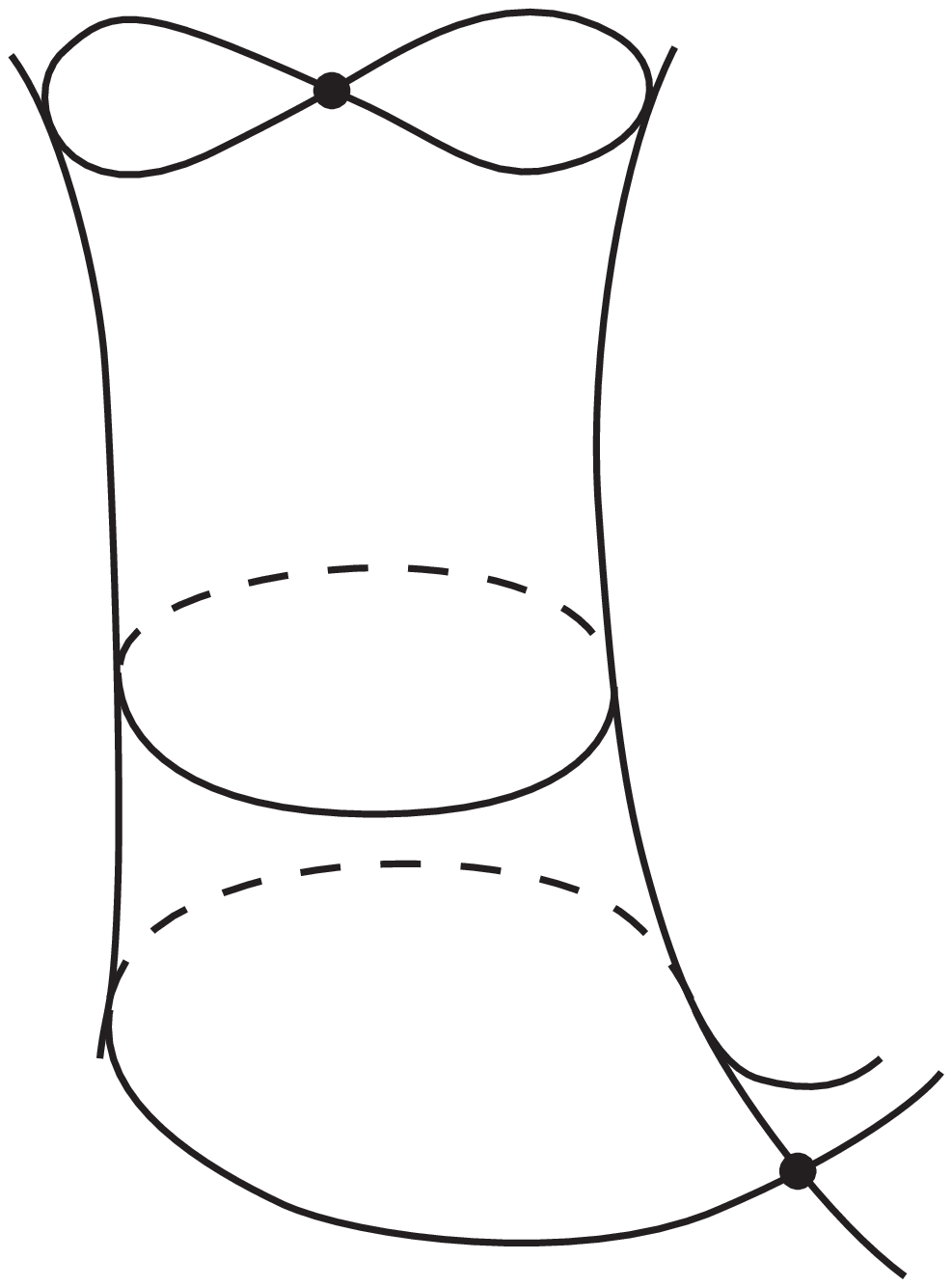,height=100pt}}
\caption{Cylinders of compact
trajectories}\label{cylinders}
\end{figure}

\begin{figure}[ht]
\centerline{\epsfig{file=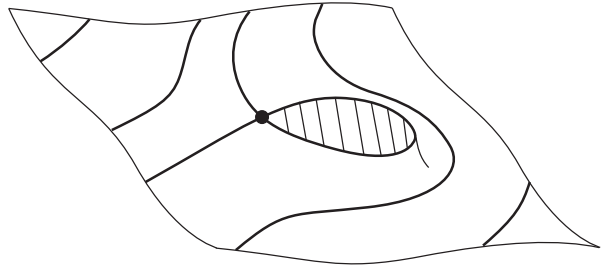,height=100pt}}
\caption{A wrapped plane}\label{wrapped}
\end{figure}

Let us describe the three-dimensional picture in more details. For
a generic level surface $M^2\subset\torus^3$, the structure of
trajectories on $M^2$ is as follows. Compact trajectories form a
few open cylinders whose ``ends'' approach either an extremum
point of the function $f|_{M^2}$ or a saddle connection cycle, see
Fig.~\ref{cylinders}.
The rest of the surface (if not empty) consists of an
even number of two-tori with or without holes, and each hole is a
saddle connection cycle. Each hole can be glued up by a planar
disk perpendicular to the vector $B$. We denote the obtained
surface by $N$. The preimage $\widehat N\subset\real^3$ of $N$
under the projection $\nu:\real^3\rightarrow\torus^3$ is a family of
finitely deformed periodic ``wrapped''  planes in $\real^3$, see
Fig.~\ref{wrapped}.

What happens to $\widehat N$ when the surface $M^2$ changes?
Small deformations of the surface $M^2$ cause just small
deformations of the tori and their covering planes. Suppose we
have a generic $1$-parametric family of surfaces $M^2(t)$. This
means that we consider a generic $1$-parametric family of
functions $f_t:\torus^3\rightarrow\real$, and for each $t$, the
surface $M^2(t)$ is defined by the equation $f_t(x)=\const$.

When the parameter $t$ varies, the connected components of $N$ are
just deformed while they stay apart from each other. However,
eventually two tori can collide and disappear or, on the contrary,
a pair of tori can be born. This occurs when $M^2(t)$ traverses a
subset which has codimension one in a natural sense. It is not
important here whether $M^2(t)$ is the family of level surface of
a single function or an arbitrary generic one-parametric family of
surfaces.

The generic tori collision was described in~\cite{D3}.
It was assumed in \cite{D3} that the family of
surfaces $M^2(t)$ is the family of level surfaces of the same
function, $M^2(t)= \{x\in\torus^3\;|\;f(x)=t\}$. However, the
argument is exactly the same for an arbitrary generic family of
surfaces.

The following two types of tori collision are possible
in the generic case.
\begin{enumerate}
\item
A cylinder of closed trajectories with bases attached to two
different components of $N$ collapses. The corresponding
codimension-one condition has the following form: two different
saddles get joined by a saddle connection.
This causes an ``interaction'' of pairs of open
trajectories lying on the collided tori, which turns
them into infinitely many closed trajectories,
see Fig.~\ref{degenerate-cylinder}.
\begin{figure}[ht]
\centerline{\epsfig{file=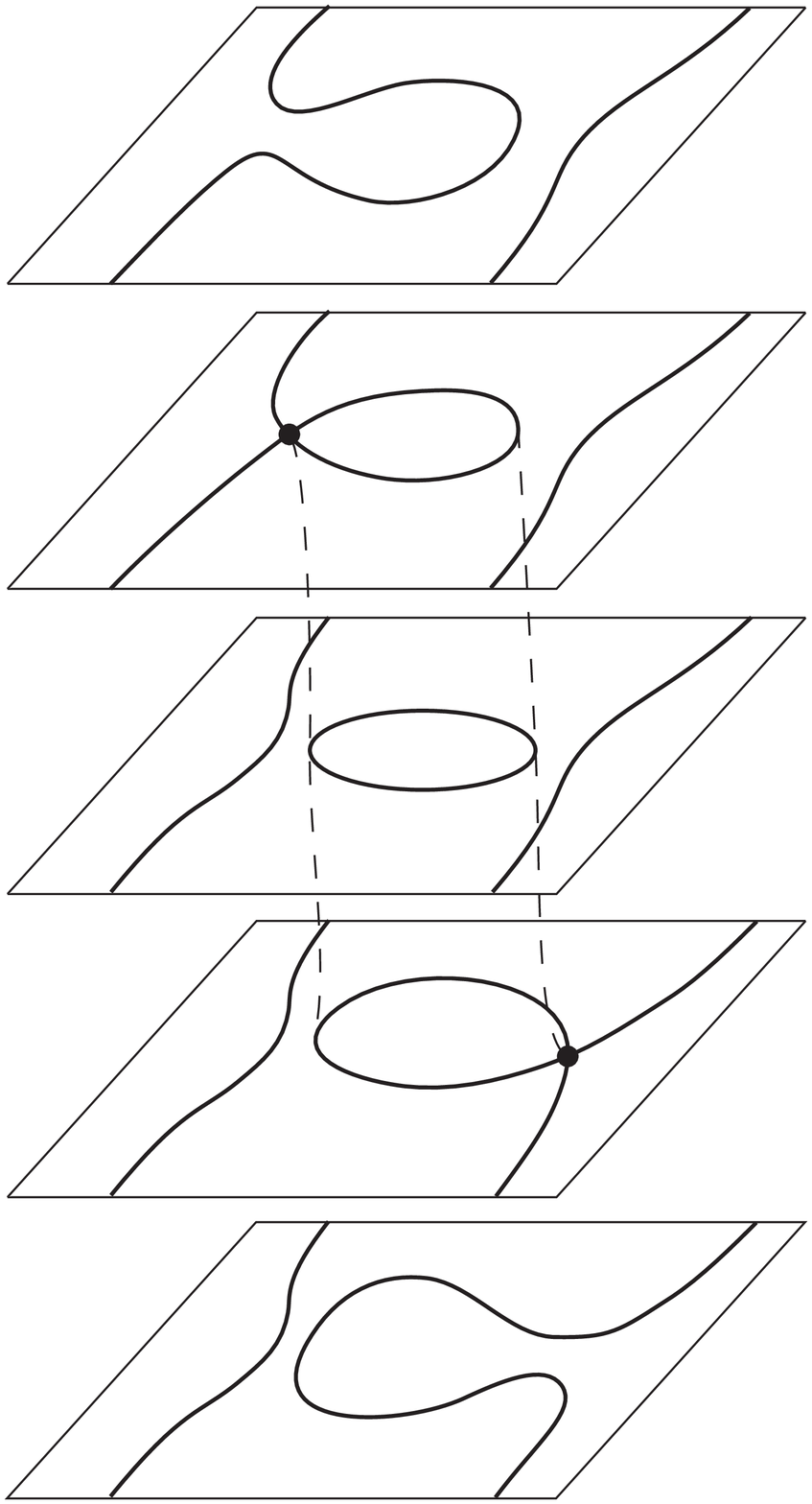,height=200pt}
\hskip0.5cm\raisebox{100pt}{$\rightarrow$}\hskip0.5cm
\epsfig{file=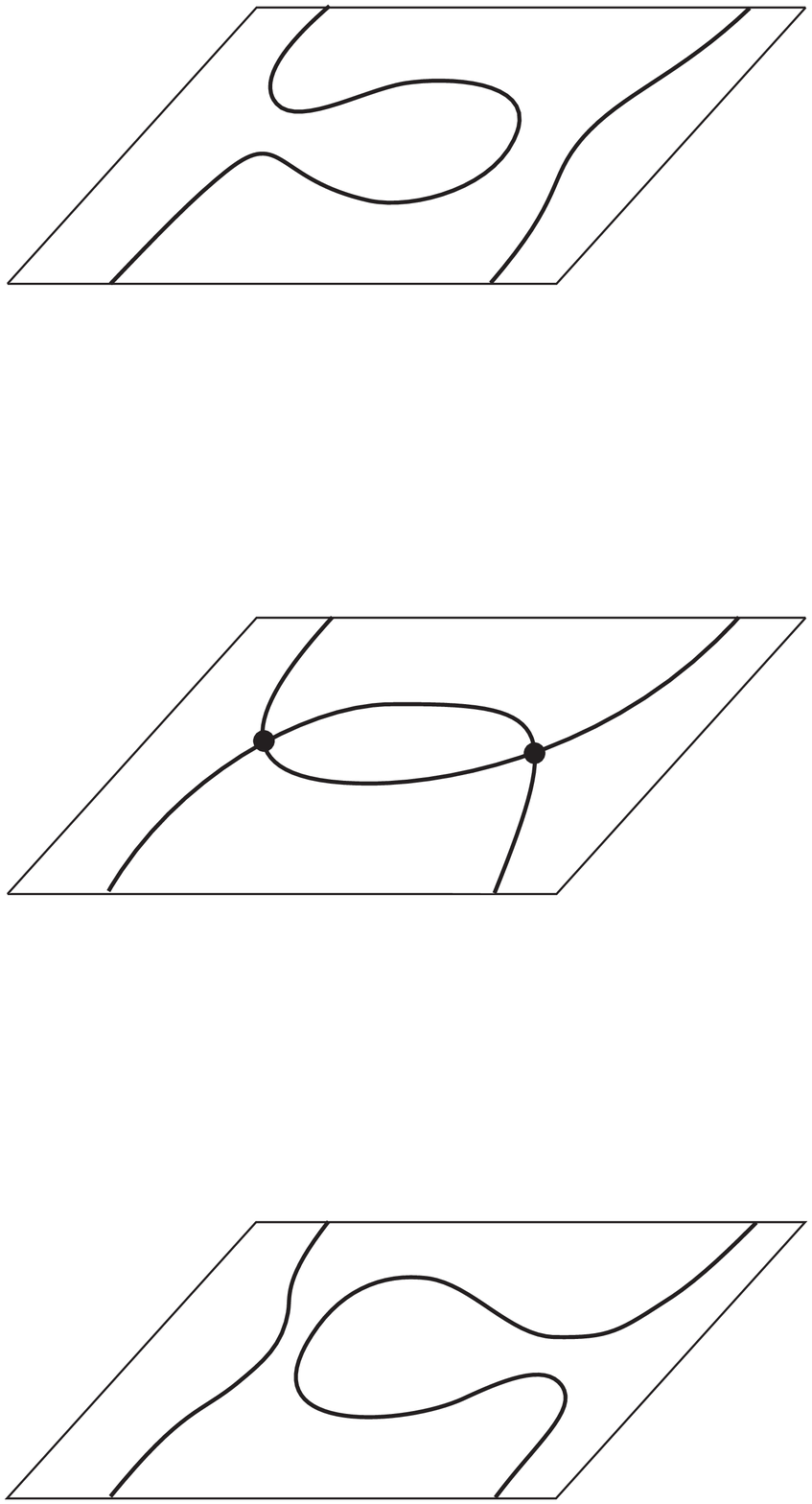,height=200pt}
\hskip0.5cm\raisebox{100pt}{$\rightarrow$}\hskip0.5cm
\epsfig{file=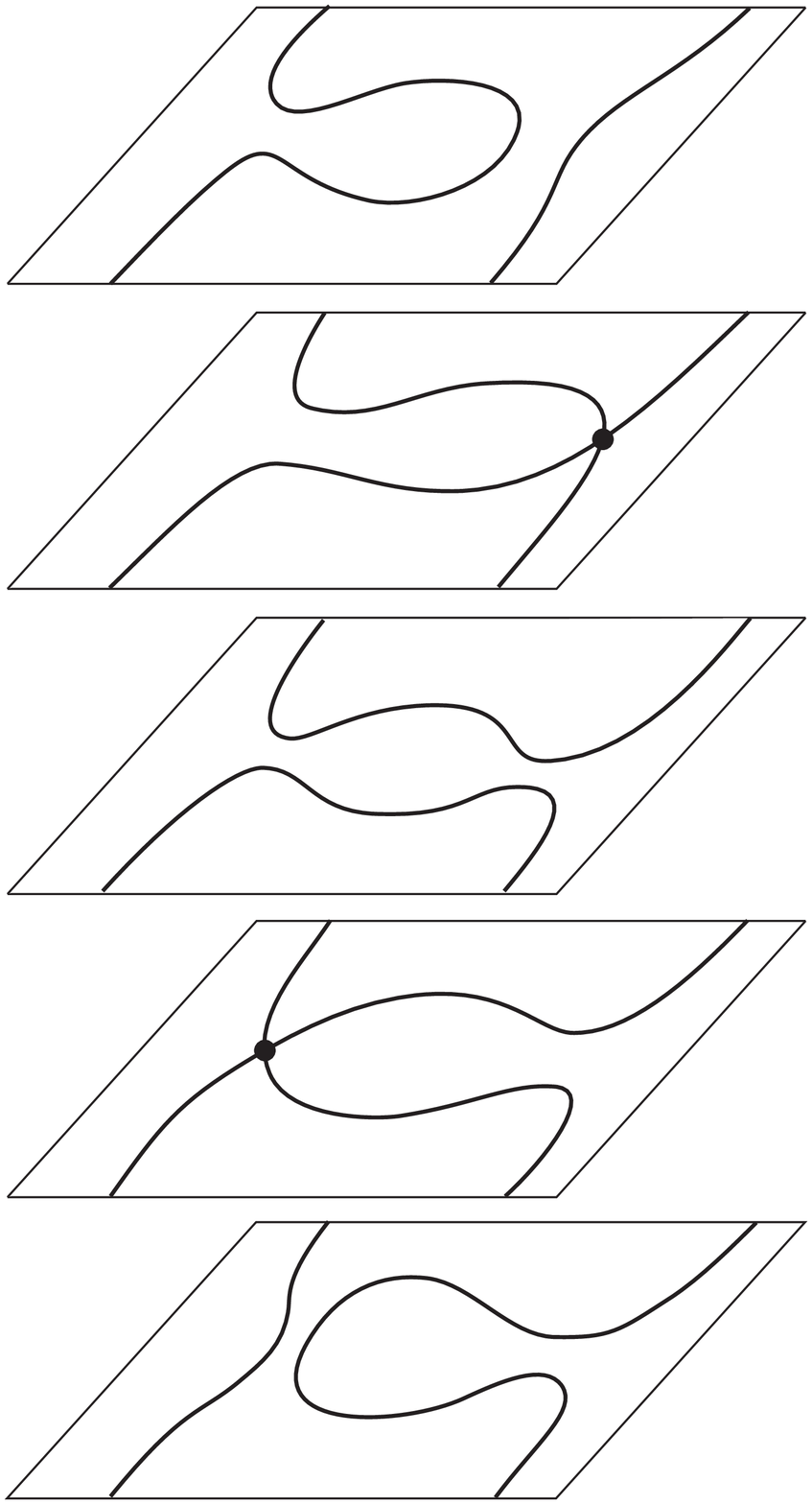,height=200pt}}
\caption{Collapse of a cylinder}
\label{degenerate-cylinder}
\end{figure}
\item
A Morse-type surgery occurs on $M^2$ that results in a one-handle added
to the surface. The behavior of trajectories is shown in Fig.~\ref{morse}.
\begin{figure}[ht]
\centerline{\epsfig{file=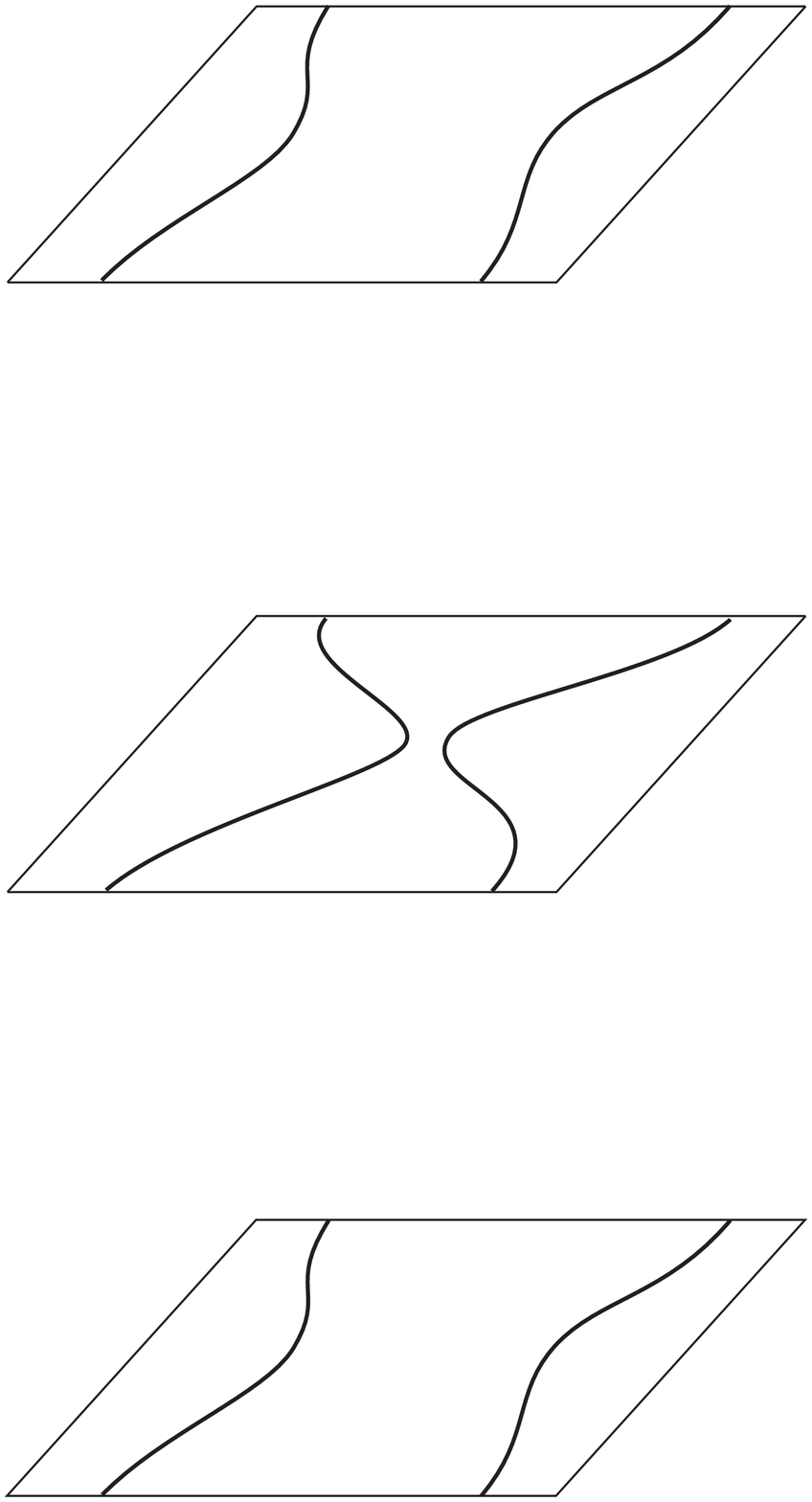,height=200pt}
\hskip0.5cm\raisebox{100pt}{$\rightarrow$}\hskip0.5cm
\epsfig{file=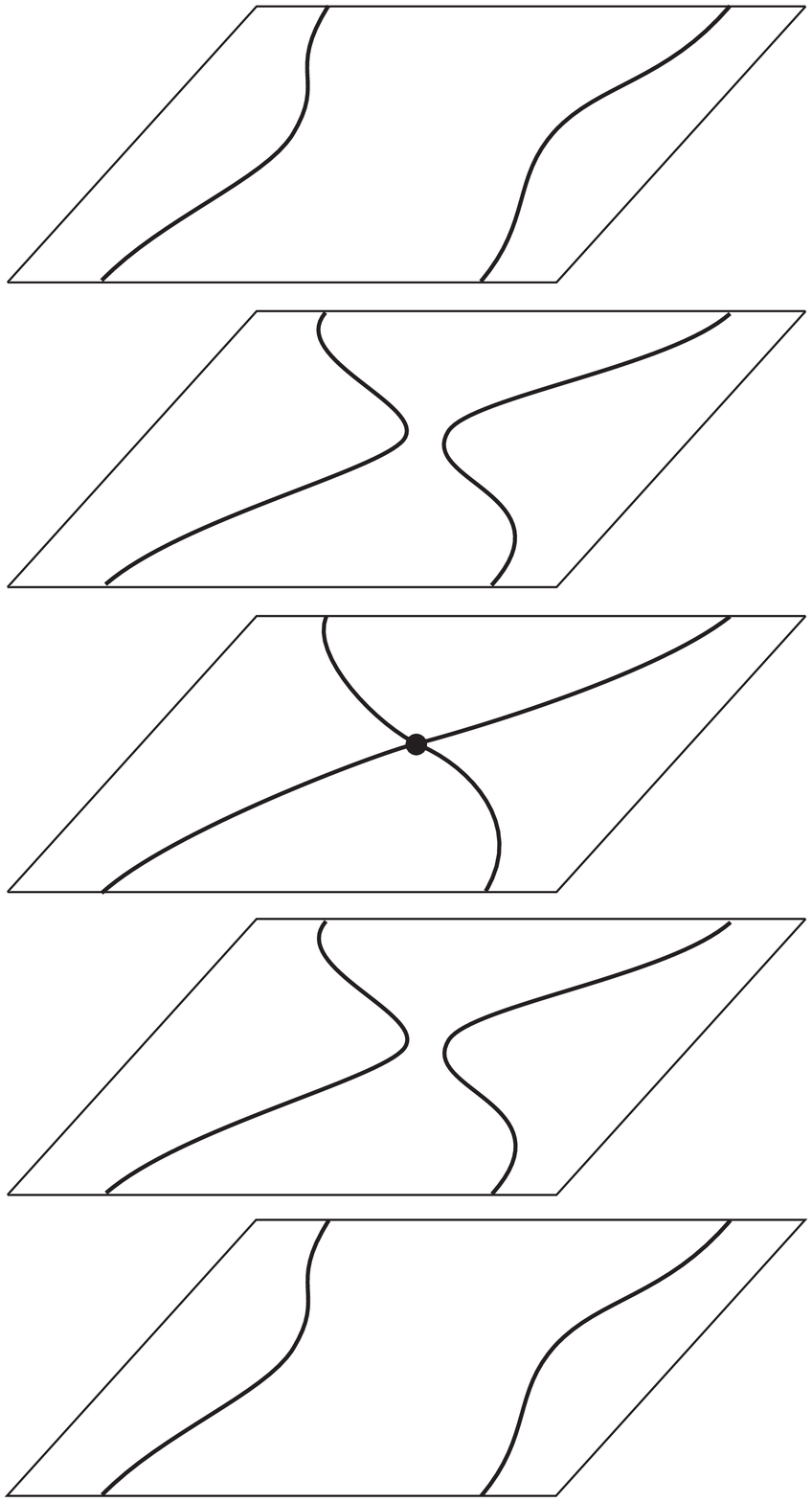,height=200pt}
\hskip0.5cm\raisebox{100pt}{$\rightarrow$}\hskip0.5cm
\epsfig{file=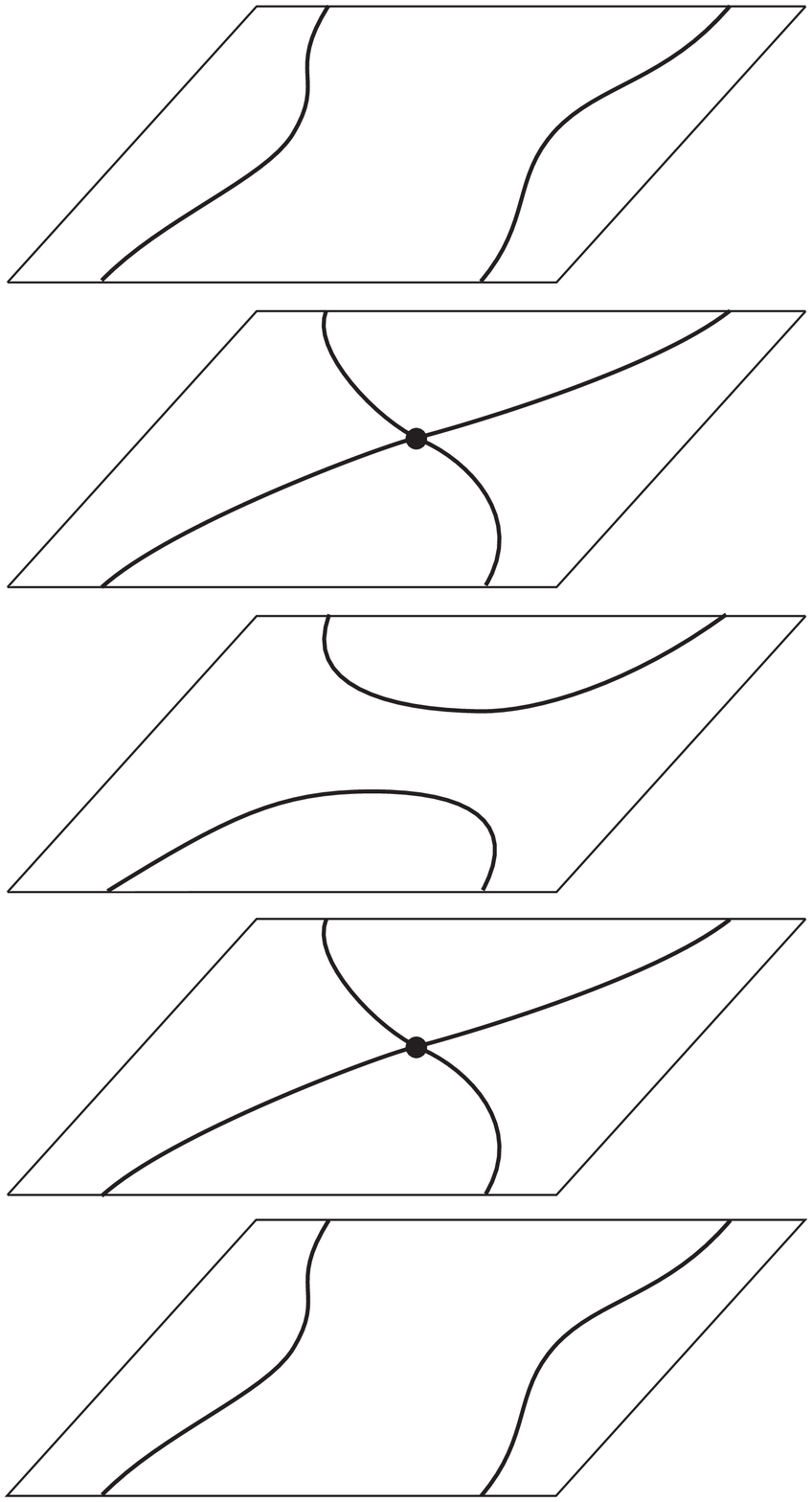,height=200pt}}
\caption{A Morse surgery destroying open trajectories}\label{morse}
\end{figure}
\end{enumerate}

It is important to note here that whenever $N$ consists of just
two tori and $M^2(t)$ passes a singularity of one of the two types
mentioned above, then all open trajectories get destroyed, so that
$N$ is empty right after the critical event. However, the
following is not proven to be impossible in a generic
one-parametric family of surfaces: bearing and canceling of a
pair of tori occurs alternatingly at moments $t_1,t_2,t_3,\dots$
so that the sequence $(t_n)$ converges to some $t_*$ and an
ergodic regime occurs on $M^2(t_*)$. The latter means that there
is an open trajectory on $M^2(t_*)$ whose closure has genus more
than one (actually, it should then be equal to three). The
existence of such ergodic regimes was proven in~\cite{D2,D3}, and it
was shown only that such regimes satisfy a codimension one
condition. However, we will not need to deal with ergodic regimes
in order to prove our result.

Now we turn to the case when $M^2(t)=M^2_t$ is the family of level
surfaces of a generic function on $\torus^3$.

Consider the restriction of the function $f$ to the plane $\real^2_{\Pi,a}$
for some $a$. Let $V\subset\real^2_{\Pi,a}$ be the union of all compact
components of $\gamma_{a,c}$ over all $c$, and $W\subset\real^2_{\Pi,a}$ the
union of all unbounded components of $\gamma_{a,c}$. We have
$\real^2_{\Pi,a}=V\cup W$, $V\cap W=\varnothing$, $V$ is open. Notice:
connected components of $V$ are not necessarily bounded. Let
$V_1,V_2,\ldots$ be the connected components of $V$. It is easy to
see that $f$ is constant on $\partial V_i$ for any $i$.

For $x\in\real^2_{\Pi,a}$ we put
$$\overline f(x)=\left\{\begin{array}{ll}f(x)&\text{if }x\in W,\\
f(\partial V_i)&\text{if }x\in V_i.\end{array}\right.$$
By doing so for all $a$, we obtain a new function $\overline f:
\real^3\rightarrow\real$. We use the following notation:
$$N_c=\{x\in\torus^3\;;\;\overline f(x)=c\}.$$
The function $f$ and its level sets $N_c$
have the following properties.

\begin{lemma}
The function $\overline f$ is a well defined continuous function
on $\torus^3=\real^3/\integer^3$.

If $U=\{c_0\}$, then $\overline f\equiv c_0$.

If $U=[c_-,c_+]$, where $c_-<c_+$, then for all but finitely many
$c\in(c_-,c_+)$, we have
$$N_c=\bigcup\limits_i\torus^2_{c,i}.$$

For all $c\in[c_-,c_+]$ a small regular neighborhood of
$N_c$ is homeomorphic to the union of a few copies of
$\torus^2\times[0,1]$.
\end{lemma}

\begin{proof}
In the case $U=\{c_0\}$ our claim is trivial.

Assume that $U=[c_-,c_+]$ with $c_-<c_+$.
By construction, $\torus^2_{c,i}\cap M_c$ consists of open trajectories,
thus, we have $\overline f(x)\equiv c$ on $\torus^2_{c,i}\cap M_c$.
The whole torus $\torus^2_{c,i}$ is obtained from
$\torus^2_{c,i}\cap M_c$ by attaching disks each of which
lies in the plane $\real^2_{\Pi,a}$ for some $a$. The boundary of such a disk is
a part of a singular unbounded component of the level set
of $f|_{\real^2_{\Pi,a}}$. By construction, we have $\overline f\equiv c$
in such a disk. Therefore, we always have
$\torus^2_{c,i}\subset N_c$.

Let us look at what happens with the tori $\torus^2_{c,i}$ when
$c$ varies. For all but finitely many $c$ the tori
$\torus^2_{c,i},\torus^2_{c,j}$ are disjoint if $i\ne j$.
Moreover, for any $c\ne c'$ and any $i,j$, the tori $\torus^2_{c,i}$
$\torus^2_{c',j}$ are always disjoint.

When $c$ varies, tori $\torus^2_{c,i}$ are continuously deformed
except at a few values of $c$, where one of the following happens:
1) two tori collide and then disappear; 2) two tori are
newly born. The latter event is opposite to the first one.

Let us describe torus collision in more detail.
At the moment of the collision we have a closed domain $W$
in $\torus^3$ that has the form of the manifold
$\torus^2\times[0,1]$ in which some intervals
$x\times[0,1]$ are collapsed to a point. There may be
just one such points $x$ or a closed disk $D^2\subset\torus^2$
of such points. The first case corresponds to an index one or index two
Morse critical point of $f$, whereas the latter
corresponds to a degenerate cylinder of closed trajectories.

The interior of the domain $W$ is filled by compact trajectories
and, by construction, the function $\overline f$ is constant
inside $W$. Thus, $W$ is a connected component of some $N_c$,
since $W$ is squeezed between the two collided tori. We call such
a $W$ \emph{pseudotorus}.

So, we have the following picture. The decomposition of $\torus^3$
into the union of (connected components of) $N_c$ over all $c$
is nothing else but a trivial fibration over $S^1$ with fibre $\torus^2$,
with a few fibres replaced by pseudotori.

Schematically this is shown in Fig.~\ref{pseudotori}.
\begin{figure}[ht]
\centerline{\epsfig{file=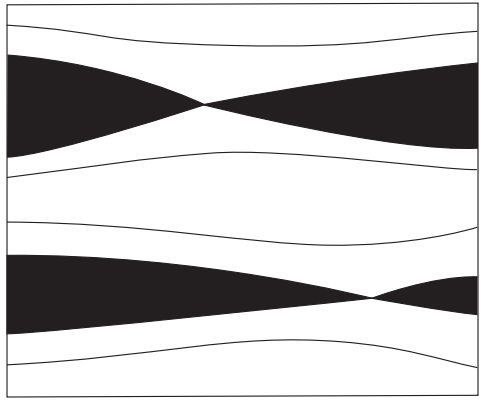,height=150pt}}
\caption{A family of tori with a few replaced by pseudotori}
\label{pseudotori}
\end{figure}

\end{proof}

Now we return to the four-dimensional case. Let $\Pi\in G_{4,2}$
be a two-plane defined by a linear mapping $B:\real^4\rightarrow\real^2$.
By $\real^2_{\Pi,a,b}$ we denote the affine plane $B^{-1}(a,b)\subset\real^4$,
and by $M^4_{\leqslant c}$ (respectively, $M^4_{\geqslant c}$)
the subset of $\torus^4$ defined by the inequality $f(x)\leqslant c$
(respectively, $f(x)\geqslant c$).

Let $N\subset\torus^4$ be a submanifold (or, more generally, a subset).
We say that $N$ is \emph{essentially below} (respectively, \emph{essentially
above}) $M^3_c$ if for any $a,b\in\real$,
the intersection $\widehat N\cap\real^2_{\Pi,a,b}$
is disjoint from all unbounded components of $\widehat M^4_{\geqslant c}
\cap\real^2_{\Pi,a,b}$ (respectively, of $\widehat M^4_{\leqslant c}\cap\real^2_{\Pi,a,b}$).
Thus, the property of $N$ to be essentially below $M^3_c$ depends on $\Pi$.

The following two facts are proved by analogy with the 3D case.

\begin{lemma}
If $N$ is essentially below or essentially
above $M^3_c$ for a given $\Pi$ then this remains
true after a small perturbation of $\Pi$, $f$, and $c$.
\end{lemma}

\begin{lemma}
If there exists a homologically nontrivial $3$-torus $N$
which is essentially above or essentially below $M^3_c$,
then the assertion of Theorem~\ref{th1} is true
for these specific $f$, $\Pi$, and $c$.
\end{lemma}

Thus, in order to prove Theorem~\ref{th1} it suffices to
show that for everywhere dense set of pairs $(f,\Pi)$,
and for each $c$,
there exists a homologically nontrivial $3$-torus $N\subset\torus^4$
which is essentially below or essentially above $M^3_c$.

Let $B=(\ell_1,\ell_2)$ be a couple of
linear functions on $\real^4$ such that
\begin{enumerate}
\item
the function $\ell_1$ is rational, i.e., $\ell_1\in(\integer^4)^*$;
\item
the restriction of $\ell_2$ to the integral three-plane $\ell_1=0$
has irrationality degree three.
\end{enumerate}
Obviously, the set of $2$-planes $\Pi=\ker B$ defined by
$\ell_1,\ell_2$ of this form is everywhere dense in $G_{4,2}$.

Without loss of generality, we assume that $\ell_1(x)=x_0$,
$\ell_2(x)=H_1x_1+H_2x_2+H_3x_3$, where
$x=(x_0,x_1,x_2,x_3)\in\real^4$, $\qopname\relax
o{rank}_\integer\langle H_1,H_2,H_3\rangle=3$. We consider the
$4$-torus $\torus^4$ as a one-parametric family of three-tori
$\torus^3_t=\{x_0=t\}$. For any $t\in[0,1]$, we deal with the
restrictions $f_t$ and $\ell_{2,t}$ of respectively $f$ and
$\ell_2$ to $\torus^3_t$ as in the three-dimensional case. We
introduce $U_t=[c_{t-},c_{t+}]$, $\overline{f_t}$, $N_{t,c}$ as
before.

Let us consider the dependence of the interval $U_t$ on $t$.
The endpoints $c_{t\pm}$ of the interval $U_t$ are continuous functions
of $t$. Moreover, in the regions where $c_{t+}>c_{t-}$
these functions are piecewise smooth. This follows from
the fact that locally, near a generic $t$, they are defined
by a condition of the form: two saddles on $M^2_{t,c_{\pm}}$
are connected by a separatrix. Here we call such intervals \emph{stability zones}.
To every stability zone there corresponds an integral vector
$\mu\in H_1(\torus^3,\integer)=\integer^3$, which we call
the \emph{label} of the zone.

\begin{figure}[p]
\noindent a)

\centerline{\begin{picture}(425,160)
\put(0,30){\epsfig{file=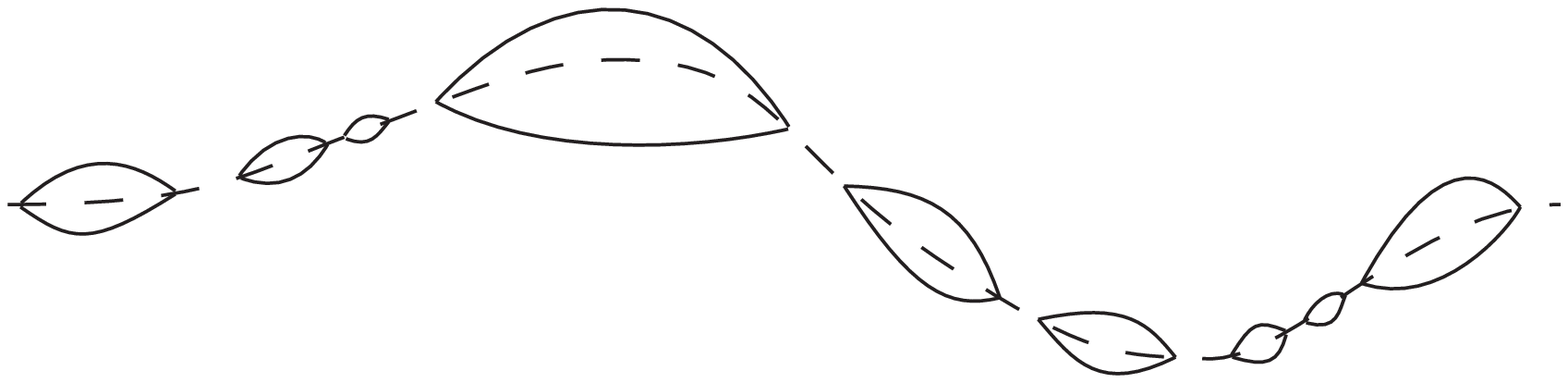,width=400pt}}
\put(0,10){\vector(1,0){425}} \put(422,0){$t$}
\put(150,80){$c_{t-}$} \put(150,130){$c_{t+}$}
\put(402,71){$c_0(t)$}
\end{picture}}

\vskip1cm
\noindent b)

\centerline{\begin{picture}(425,160)
\put(0,30){\epsfig{file=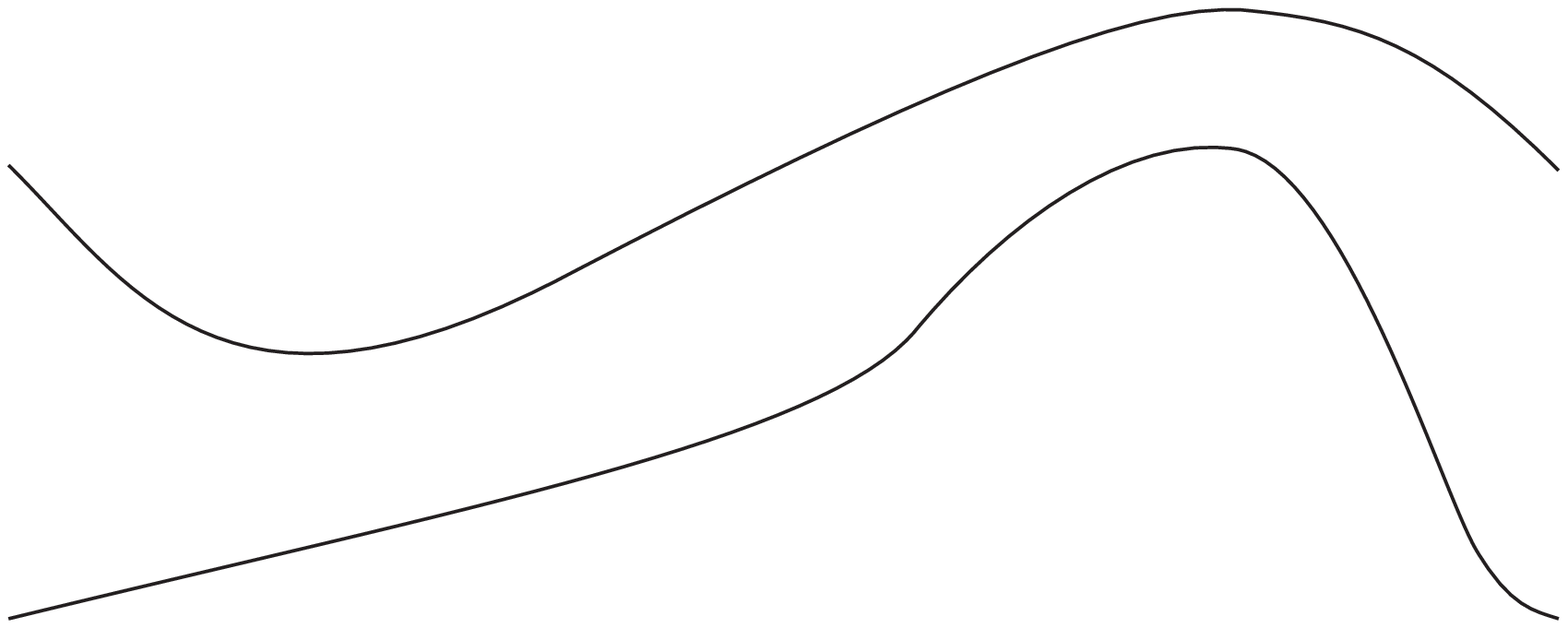,width=400pt}}
\put(0,10){\vector(1,0){425}} \put(422,0){$t$}
\put(150,76){$c_{t-}$} \put(150,135){$c_{t+}$}
\multiput(0,152)(10,0){41}{\line(1,0)5}
\put(408,149){$(c_-)_{\mathrm{max}}$}
\multiput(0,99)(10,0){41}{\line(1,0)5}
\put(408,96){$(c_+)_{\mathrm{min}}$}
\end{picture}}

\vskip1cm
\noindent c)

\centerline{\begin{picture}(425,160)
\put(0,30){\epsfig{file=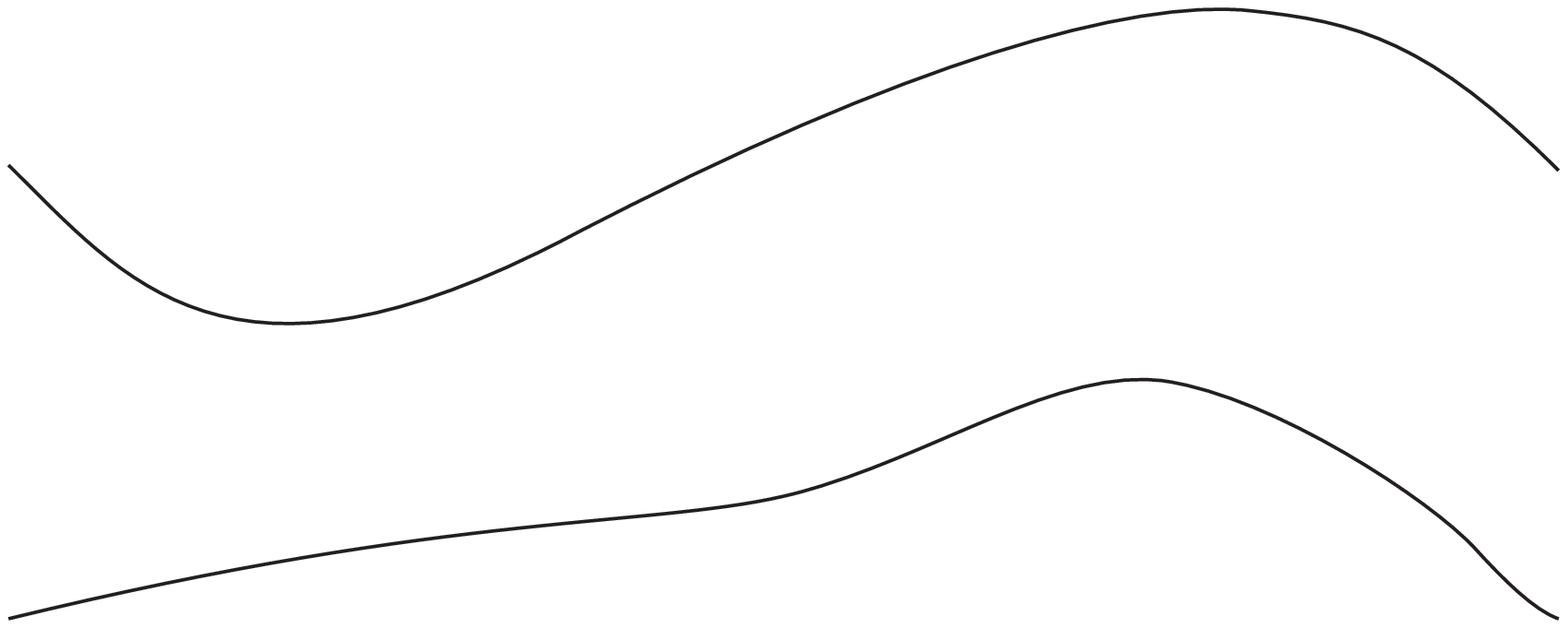,width=400pt}}
\put(0,10){\vector(1,0){425}} \put(422,0){$t$}
\put(150,62){$c_{t-}$} \put(150,145){$c_{t+}$}
\multiput(0,106)(10,0){41}{\line(1,0)5}
\put(408,89){$(c_-)_{\mathrm{max}}$}
\multiput(0,92)(10,0){41}{\line(1,0)5}
\put(408,103){$(c_+)_{\mathrm{min}}$}
\end{picture}}
\caption{Functions $c_{t\pm}$ in the generic case}\label{cpm}
\end{figure}

Figure~\ref{cpm} shows how the functions $c_{t\pm}$ may look like
in the generic case.
It is possible that at some $t$ we have $c_{t+}=c_{t-}$,
see Fig.~\ref{cpm}a). This may
occur at the boundary of a stability zone or at $t$ such that
the open trajectories in $M^3_t$ are chaotic.
For all such $t$ we have $c_{t+}=c_{t-}=c_0(t)$, where $c_0(t)$ is
a piecewise smooth function of $t$. It is defined locally by a
condition of the form: the sum of ``heights'' of certain saddles
equals to zero, see~\cite{D2,D3}.

It is most likely that any ``chaotic" $t$ must be an accumulating
point of an infinite sequence of stability zones. In other words,
it cannot happen that the equality $c_{t+}=c_{t-}$ holds
everywhere in a nontrivial interval $(t_1,t_2)$. However, this
does not follow directly from the previous works~\cite{D2,D3}, and
will not be used here.

\begin{lemma}\label{=/=}
The equality
\begin{equation}\label{c+=c-}
\min_tc_{t+}=\max_tc_{t-}
\end{equation}
does not hold for a generic $f$.
\end{lemma}

\begin{proof}
We should consider the following three cases.

Case 1. For all $t$ we have $c_{t+}>c_{t-}$. Then there exists
a smooth periodic function $g(t)$ such that
$c_{t+}>g(t)>c_{t-}$. Condition~(\ref{c+=c-})
will not hold if we disturb the function $f$ in the following way:
$f(t,x_1,x_2,x_3)\mapsto f_\varepsilon(t,x_1,x_2,x_3)=f(t,x_1,x_2,x_3)+
\varepsilon g(t)$, where $|\varepsilon|$ is sufficiently small.
So,~(\ref{c+=c-}) impose a codimension one condition on $f$ in this case.

Case 2. For two different $t=t_1,t_2$ we have $c_{t+}=c_{t-}$.
Then condition~(\ref{c+=c-}) will not hold after an arbitrary
perturbation $f\mapsto f+\varepsilon g(t)$, where $g(t)$
is an arbitrary function with $g(t_1)\ne g(t_2)$.

Case 3. There is exactly one $t=t_0$ such that $c_{t+}=c_{t-}=c_0$ holds,
and we have $c_{t+}>c>c_{t-}$ for all $t\ne t_0$.
Let us take $t$ close to $t_0$. The interval $[c_{t-},c_{t+}]$
is then small, which means that, when $c$ varies from
$(c_{t-}-\delta)$ to $(c_{t+}+\delta)$ with $\delta>0$,
a pair of tories $N_{t,c}$ is born at $c=c_{t-}$ and then almost
immediately destroyed at $c=c_{t+}$. Therefore,
there are two cylinders of closed trajectories on $M_{t,c_0}$
of very small height.

Now let us fix $c=c_0$ and vary $t$. When $t$ approaches $t_0$
from the left, say, we will have two closed trajectory cylinders
that get degenerate at the moment $t=t_0$. When $t$ passes $t_0$,
two cylinders must appear again. The main point here is that
those, new, cylinders appear from \emph{the same} pair of degenerate
cylinders. Indeed, the pair of degenerate cylinders that we obtain
when $t$ approaches $t_0$ from the left cuts $M_{c_0,t_0}$
into two tori. Since the irrationality degree of $\ell_2$
is equal to three, one can show that there may no other degenerate cylinder on
$M_{t_0,c_0}$.

Thus, we have the following picture. When $t$ goes from $t_0-\delta$
to $t_0+\delta$, two closed trajectory cylinders degenerate
and then regenerate again. Let $h_1(t)$, $h_2(t)$ be there heights.
So, we not only have $h_1(t_0)=h_2(t_0)=0$, but also
$h_1'(t_0)=h_2'(t_0)=0$, which impose a codimension two
condition on the function $f$.
\end{proof}

By construction we have

\begin{lemma}\label{<>}
For any $t$ and $c>c_{t+}$ (respectively, $c<c_{t-}$),
the torus $\torus^3_t$ is essentially
below (respectively, essentially above) $M^3_c$.
\end{lemma}

According to Lemma~\ref{=/=} only the following two cases are possible:
1) $\min_tc_{t+}<\max_tc_{t-}$; 2) $\min_tc_{t+}>\max_tc_{t-}$.
In Case~1, for any $c$ we have either $c>\min_tc_{t+}$
or $c<\max_tc_{t-}$, and by Lemma~\ref{<>} we are done.

So, it remains to consider Case~2,
$\min_tc_{t+}>\max_tc_{t_-}$. This inequality means,
in particular, that the intervals $U_t$
have a nontrivial intersection $U=\cap_tU_t=[c_-,c_+]$, and we
have just one stability zone, which covers the whole circle $S^1$.
This is illustrated in Fig.~\ref{cpm}c).

We define a function $\overline f:\torus^4\rightarrow\real$
as in the three-dimensional case by considering
intersections $M_c^3\cap\real^2_{\Pi,a,b}$, and introduce notation
$$N_c=\{x\in\torus^4\;;\;\overline f(x)=c\}.$$
By construction, $\overline{f_t}$
coincides with the restriction of $\overline f$ to $\torus^3_t$, and we
have
$$N_c=\bigcup\limits_tN_{t,c}.$$

For almost
any $c$, $t$, the intersection $N_c$ with $\torus^3_t$ is either empty
or consists
of $2$-tori, and all those tori have the same up to sign homology class
$\alpha\in H_2(\torus^3,\integer)$.
In this case, the whole torus $\torus^4$ has the structure
of a trivial $\torus^2$-bundle over $\torus^2$ with a $1$-parametric
family of fibres replaced by pseudotori.
\begin{figure}[ht]
\centerline{\begin{picture}(400,320)
\put(0,20){\epsfig{file=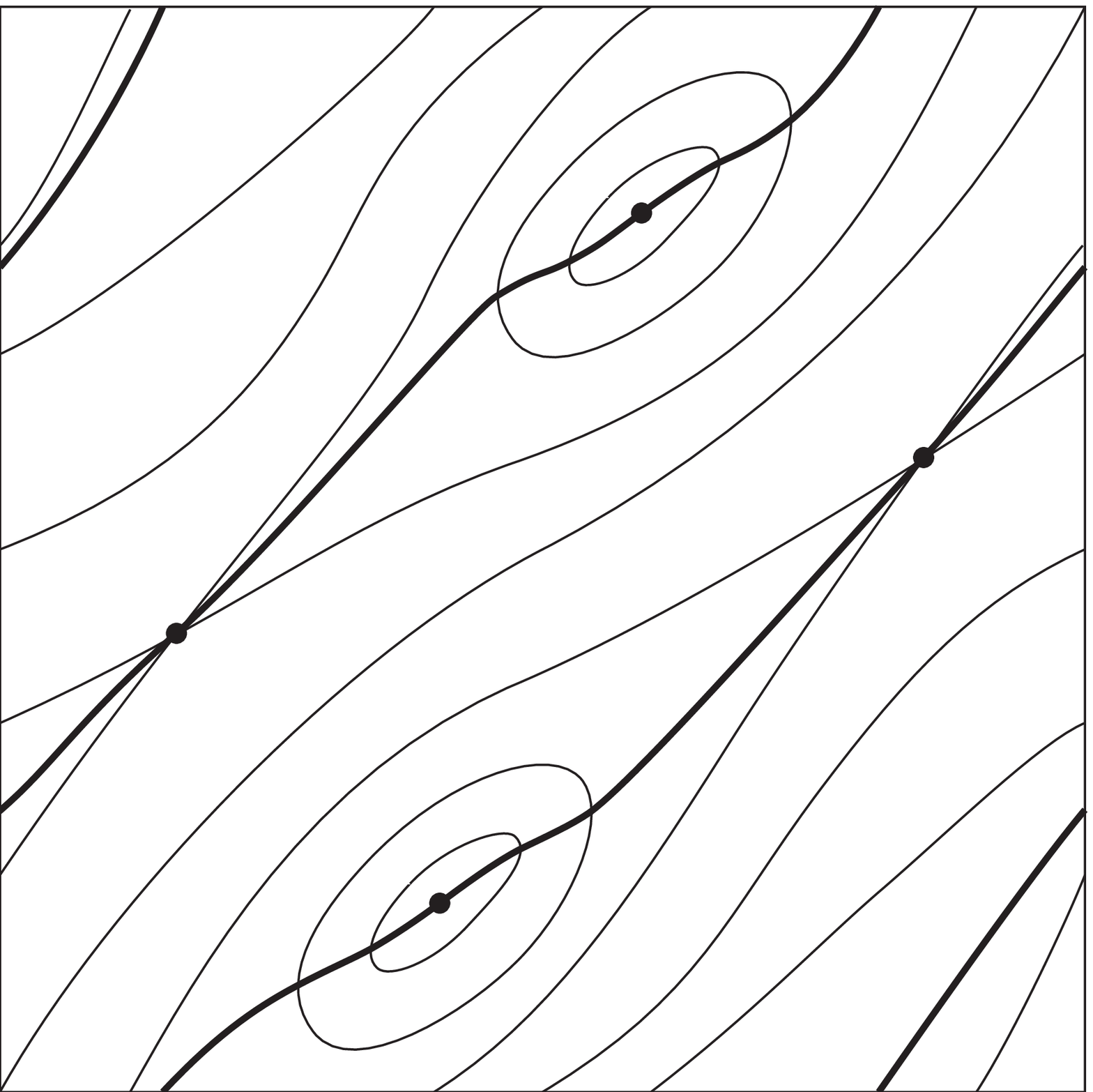,height=300pt}}
\put(0,10){\vector(1,0){310}} \put(305,0){$t$} \put(400,160){\hbox
to0pt{\hss Pseudotori}} \put(340,170){\vector(-1,1){55}}
\put(340,160){\vector(-2,-3){37}}
\end{picture}}
\caption{Level lines of $\overline f$ on $\torus^2$}\label{t2}
\end{figure}
The function
$\overline f$ is constant over each fibre, so it can be considered
as a function on $\torus^2$.

Figure~\ref{t2} illustrates the structure of level lines of $\overline f$
viewed as a function on $\torus^2$. The preimage of a generic point
is a $2$-torus imbedded in $\torus^4$, and all those $2$-tori are
``parallel". Whenever the function $\overline f$ has an extremum on
the line $t=\mathrm{const}$, the preimage of the critical point
is a pseudotorus.

Including a pseudotorus into a $1$-parametric family of $2$-tori
does not change the topological type of the union of the tori. Indeed,
as we mentioned above, a small regular neighborhood of a pseudotorus
in $\torus^3$ is homeomorphic to $\torus^2\times(0,1)$. In other words,
attaching collars $\torus^2\times(0,1)$ to a pseudotorus again
gives $\torus^2\times(0,1)$.

We conclude the following from this.

\begin{lemma}
For almost all $c$ each connected component
of $N_c$ will be homeomorphic to $\torus^3$. The exceptions are those
$c$ that are critical values of the function $\overline f$ on $\torus^2$.
\end{lemma}

For a generic $f$ we obtain a generic Morse function $\overline f$
on $\torus^2$. For such a function, there must be an interval $[c_1,c_2]$
such that, whenever we have $c\in[c_1,c_2]$,
the level line $\overline f=c$ contains a closed curve non-homologous to
zero in $\torus^2$. The preimage $N_c$ of this level line in $\torus^4$
is a $3$-torus non-homologous to zero.
Thus, we get the following.

\begin{lemma}
There exist $c_1$, $c_2$ such that $c_1<c_2$ and both $N_{c_1}$
and $N_{c_2}$ contain a connected component homeomorphic to $\torus^3$
and non-homologous to zero.
\end{lemma}

It remains to notice that whenever $c>c_1$ the hypersurface $N_{c_1}$
is essentially below $M^3_c$, and whenever $c<c_2$ the hypersurface
$N_{c_2}$ is essentially above $M^3_c$. Thus, for all $c$
we have a non-homologous to zero $3$-torus which is either essentially
above or essentially below $M^3_c$, and we are done in Case~2.

\end{document}